\documentclass[12pt]{article}

\usepackage[all]{xy}
\usepackage[dvips]{graphicx}
\usepackage{enumerate,hyperref}
\usepackage{ifthen,color}
\usepackage{amsmath,amsfonts,amssymb,amsthm,dsfont}
\usepackage{mathrsfs}

\newtheorem{thm}{Theorem}[section]
\newtheorem{prop}[thm]{Proposition}
\newtheorem{coro}[thm]{Corollary}

\newtheorem{exa}[thm]{Example}
\newtheorem{lem}[thm]{Lemma}

\newtheorem{de}[thm]{Definition}

\newcommand{\comp}{ \,{\scriptstyle \stackrel{\circ}{}}\, }

\newcommand{\Acal}{\mathcal{A}}

\newcommand{\Bcal}{\mathcal{B}}

\newcommand{\Hcal}{\mathcal{H}}
\newcommand{\calH}{\mathcal{H}}

\newcommand{\RN}{ \mathbb{R} } 
 
\newcommand{\CN}{ \mathbb{C} } 
 
\newcommand{\HN}{\mathbb{H}}
\newcommand{\NN}{\mathbb{N}}

\newcommand{\eps}{\epsilon}

\newcommand{\mapmap}{\textbf{map}}

\newcommand{\MatA}{\text{M}_{\Acal}}

\newcommand{\EndA}{\mathcal{R}_{\Acal}}

\newcommand{\MM}{\mathbf{M}}
\newcommand{\Acalx}{\mathcal{A}^{\times}}
\newcommand{\Acalzd}{\mathbf{zd}(\mathcal{A})}

\newcommand{\ds}{\displaystyle}

\newcommand{\Abound}{m_{\Acal}}

\begin{document}
\begin{titlepage}

\title{{\sc Theory of Series in the $\Acal$-calculus and the $N$-Pythagorean Theorem}}
\author{Daniel Freese\space\space\space\space\space\space\space James S. Cook \\ djfreese@iu.edu \space\space \space jcook4@liberty.edu}


\maketitle

\begin{abstract}
In this paper we study sequences, series, power series and uniform convergence in the $\mathcal{A}$-Calculus. Here $\mathcal{A}$ denotes an associative unital real algebra. We say a function is $\mathcal{A}$-differentiable if it is real differentiable and its differential is in the regular representation of the algebra. We show the theory of sequences and numerical series resembles the usual theory, but, the proof to establish this claim requires modification of the standard arguments due to the submultiplicativity of the norm on $\mathcal{A}$. In contrast, the theorems concerning divergence of power series over $\mathcal{A}$ are modified notably from the standard theory. We study how the ratio, root and geometric series results are modified due to both the submultiplicativity of the norm and the calculational novelty of zero-divisors. We establish the usual calculations with power series transfer nicely to the $\mathcal{A}$-calculus. Power series are used to define sine, cosine, hyperbolic sine, hyperbolic cosine and the exponential. Finally, special functions are introduced and we derive the $N$-Pythagorean Theorem.
\end{abstract}

\vspace{0.2in}

\end{titlepage}

\section{Introduction and overview}
We use $\Acal$ to denote a real unital associative algebra of finite dimension. Elements of $\Acal$ are known as $\Acal$-numbers. We study calculus where real numbers have been replaced by $\Acal$-numbers. The resulting calculus we refer to as $\Acal$-calculus. Our typical goal is to find theorems which apply to as large a class of real associative algebras as possible.  To our knowledge, the results presented in Sections 3 to 7 have not appeared in the literature in the generality which we supply in this work. \\

\noindent
Setting aside complex analysis, calculus over more general number systems have been studied from about 1890 to the present time. There are too many papers to list. To see how our current framework relates to the existing literature please see \cite{cookAcalculusI}. \\

\noindent
The goal of this paper is to study sequences, series and power series over an algebra. This serves as a foundation for an ongoing project to develop $\Acal$-calculus generalizations of the standard calculational techniques. We should mention N. BeDell has written three supplementary algebra papers \cite{bedellI},\cite{bedellII} and \cite{bedellIII} which provide further algebraic discussion of zero-divisors, logarithms and the $N$-Pythagorean Theorem proved in Section \ref{sec:nthagtheorem} of this paper. Then in the sequel to this paper \cite{cookbedell} one of the authors and N. BeDell present the theory of $\Acal$-ordinary differential equations. \\

\noindent
In Section \ref{sec:reviewofAcalculus} we review the major developments from \cite{cookAcalculusI}. In particular, we discuss representations, submultiplicativity of the norm, the definition and theory of differentiation over an algebra and select theorems from integral calculus over $\Acal$. Many further examples, proofs and motivations can be found in \cite{cookAcalculusI} and we encourage the reader to consult that paper before digesting this current work. \\

\noindent
We begin the presentation of new results in Section \ref{sec:sequences} where we develop the theory of sequences over $\Acal$ following \cite{Rudin}. We show the usual arithmetic of limits transfers nicely to sequences in $\Acal$. The results are quite natural, however, the submultiplicativity complicated the usual proofs from real or complex analysis. \\

\noindent
Numerical series over $\Acal$ are covered in Section \ref{sec:series}. We found that the usual elementary convergence and divergence tests are meaningful over an algebra. The $n$-th term test, comparison test, absolute convergence, root and ratio test all naturally generalize over $\Acal$. The Cauchy Criterion is meaningful and the Cauchy product exists to multiply series where at least one is absolutely convergent. Once more, the proofs required significant modification due to the submultiplicativity of the norm.  \\

\noindent
In Section \ref{sec:powerseries} we study power series in $\Acal$. We attempt to generalize the elementary convergence and divergence theorems of real calculus. We find the Root and Ratio Test need significant modification. One aspect of the modification is that the submultiplicative constant $\Abound$ appears in convergence results. For example, where $R= 1/ \alpha$ in the usual calculus we found $R = 1/ \Abound \alpha$. Also, the geometric series $1+z + z^2+ \cdots $ converges for $\| z \| < 1/ \Abound$. A second, and initially perplexing, modification is seen in the absense of divergence cases for the Root and Ratio Tests for power series over $\Acal$. The multiplicative property of the norm in real or complex analysis is important to obtain the boundary between divergence and convergence of a power series. Submultiplicativity of the norm on $\Acal$ spoils the usual argument for divergence for both the Root and Ratio Tests. However, we also understand the reason for this modifcation in view of the phenomenon seen in Example \ref{exa:bandconverge}. Zero divisors allow new domains of convergence which are not seen in single-variate analysis over $\RN$ or $\CN$. Uniform convergence and Weierstrauss $M$-Test are studied. The standard theorems concerning sequences of uniformly convergent functions hold for $\Acal$-differentiable sequences of functions.  The integral of the limiting function is the limit of the integrals, however, the result for derivatives is not as simple. See Theorem \ref{thm:difflimitfunseries} which mirrors the usual theorem of real analysis. We show the term-by-term derivative of a power series in $\Acal$ is indeed the derivative of the given power series. Given an entire function on $\RN$ we show there exists a unique entire extension to $\Acal$. We show entire functions are absolutely convergent on $\Acal$ and uniformly convergent on any finite ball in $\Acal$. Finally, we show the product of entire functions is again entire. Indeed, the entire functions on $\Acal$ form an algebra. \\

\noindent
Transcendental functions such as the exponential, sine, cosine and hyperbolic sine and cosine are covered in Section \ref{sec:transfunct}. Theorem \ref{thm:entire} indicates the defintions we offer are inescapable. We find the usual indentities for the elementary functions on $\RN$ extend naturally to any commutative unital algebra $\Acal$. \\

\noindent
Section \ref{sec:nthagtheorem} reverses the direction of study from that of Section \ref{sec:transfunct}. We ask, given a specific choice of $\Acal$, which functions appear naturally? In particular, we study functions which appear as component functions of the exponential. We call these the {\it special functions} of $\Acal$. We find a theorem we call the $N$-Pythagorean Theorem which provides an identity which holds for the special functions. This theorem makes the identities $\cos^2 \theta + \sin^2 \theta =1$ and $\cosh^2 \phi - \sinh^2 \phi = 1$ part of a chain of such identities.

\section{Review of $\Acal$ calculus} \label{sec:reviewofAcalculus}
\noindent
In this Section we have two main goals. First, to provide necessary background to understand the new theory developed in the later sections. Second, to alert the reader to some of the major results which are already established in  \cite{cookAcalculusI}. Please consult \cite{cookAcalculusI} for references and discussion of how our work connects to the existing literature. 

\subsection{Algebra and the regular representations}
We say\footnote{we write $(\Acal, \star)$ to emphasize the pairing where helpful} $\Acal$ is an {\bf algebra} if $\Acal$ is a finite-dimensional real vector space paired with a function $\star: \Acal \times \Acal \rightarrow \Acal$ which is called {\bf multiplication}. In particular, the multiplication map satisfies the properties below: 
\begin{enumerate}[{\bf (i.)}]
\item {\bf bilinear:} $ (cx + y ) \star z = c( x \star z)+ y \star z$ and $x \star (cy+ z) = c(x \star y)+ x \star z $ for all $x,y,z \in \Acal$ and $c \in \RN$,
\item {\bf associative:} for which $x \star (y  \star z) = (x \star y)  \star z$  for all $x,y,z \in \Acal$ and,
\item {\bf unital:} there exists $\mathds{1} \in \Acal$ for which $\mathds{1} \star x = x$ and $x \star \mathds{1} = x$.
\end{enumerate}
We say $x \in \Acal$ is an {\bf $\Acal$-number}. If $x \star y = y \star x$ for all $x,y \in \Acal$ then $\Acal$ is {\bf commutative}. \\

\noindent
The {\bf left-multiplication by $x$} is the map $L_x: \Acal \rightarrow \Acal$ defined by $L_x(y) = x \star y$ for all $y \in \Acal$. Observe, by associativity of $\Acal$,
\begin{equation}
 L_x(y) = x \star 1 \star y = L_x(1) \star y.
\end{equation}

\noindent
A linear transformation $T: \Acal \rightarrow \Acal$ is {\bf right $\Acal$ linear} if $T(x \star y) = T(x) \star y$ for all $x, y \in \Acal$. We say the set $\EndA$ of all right $\Acal$ linear transformations forms the {\bf regular representation} of $\Acal$. Since $\Acal$ is unital the regular representation is isomorphic\footnote{isomorphic as associative real algebras, we say $(\Acal, \star)$ and $(\Bcal, \comp)$ are isomorphic if there is a linear bijection $\Psi: \Acal \rightarrow \Bcal$ for which $\Psi( x \star y) = \Psi(x) \comp \Psi(y)$ for all $x,y \in \Acal$.}to $\Acal$. The isomorphism from $\Acal$ to $\EndA$ is denoted $\mapmap$ and we find it convenient to use $\#$ for $\mapmap^{-1}$. In particular, 
\begin{equation}
\mapmap(x) = L_x \qquad \& \qquad \#(T) = T(1).
\end{equation}
The idea here is that $\#(T)$ provides the $\Acal$ number which corresponds to $T$.
If $\beta$ is a basis for $\Acal$ then the {\bf matrix regular representation} of $\Acal$ with respect to $\beta$ is
\begin{equation}
 \MatA(\beta) = \{ [T]_{\beta, \beta} \ | \ T \in \EndA \}
\end{equation} 
where $[T]_{\beta,\beta}$ denotes the matrix of $T$ with respect to the basis $\beta$. In the case $\Acal = \RN^n$ we may forego the $\beta$ notation and write
\begin{equation}
 \MatA = \{ [T] \ | \ T \in \EndA \}
 \end{equation}
for the {\bf regular representation} of $\Acal$. There is a natural isomorphism of $\Acal$ and $\MatA$: If $\beta = \{ v_1 , \dots , v_n \}$ is a basis for $\Acal$ where $v_1 = \mathds{1}$ then 
\begin{equation}
\MM( x) = \left[ [x]_{\beta}| [x \star v_2]_{\beta}| \cdots | 
[x \star v_n]_{\beta} \right]
\end{equation}
where $[x]_{\beta}$ is the coordinate vector of $x$ with respect to $\beta$. In many applications we consider the case $\Acal = \RN^n$ with $\beta = \{ e_1, \dots , e_n \}$ the usual standard basis such that $e_1 = \mathds{1}$. Given these special choices we obtain much improved formula
\begin{equation}
 \MM(x) = [x| x \star e_2| \cdots | x \star e_n]. \end{equation}

\begin{exa} \label{Ex:number1}
The {\bf complex numbers} are defined by $\CN = \RN \oplus i\RN$ where $i^2=-1$. We denote $a+ib = [a,b]^T$ corresponding to our identifications $e_1=1$ and $e_2=i$. Notice, $(a+ib)e_2 = (a+ib)i =ia-b = [-b,a]^T$. Thus, $\mathbf{M}(a+ib) = \left[ \begin{array}{cc} a  & -b  \\  b & a \end{array} \right]$ is a typical element of the regular matrix representation of $\CN$ which we denote $\text{M}_{\CN}$.
\end{exa}

\begin{exa} \label{Ex:number2}
The {\bf hyperbolic} numbers are given by $\calH = \RN \oplus j\RN$ where $j^2=1$. Identifying $e_1=1$ and $e_2=j$ we have $a+bj = [a,b]^T$. Moreover, 
\begin{equation}
(a+bj)e_2 = (a+bj)j = aj+b = [b,a]^T. 
\end{equation}
Therefore, $ \mathbf{M}(a+bj) = \left[ \begin{array}{cc} a  & b  \\  b & a \end{array} \right]$ is a typical matrix in $\text{M}_{\calH}$.
\end{exa}

\noindent
We say $x \in \Acal$ is a {\bf unit} if there exists $y \in \Acal$ for which $x \star y = y \star x = \mathds{1}$. The set of all units is known as the {\bf group of units} and we denote this by $\Acalx$. We say $a \in \Acal$ is a {\bf zero-divisor} if $a \neq 0$ and there exists $b \neq 0$ for which $a \star b = 0$ or $b \star a =0$. Let $\Acalzd = \{ x \in \Acal \ | \ x=0 \ \text{or $x$ is a zero-divisor} \}$.

\begin{exa} \label{prop:isomorphismhyperbolic}
 and $\mathbf{zd}( \calH ) = \{ a+bj \ | \ a^2=b^2 \}$ whereas $\calH^{\times} = \{ a+bj \ | \ a^2 \neq b^2 \}$. The reciprocal of an element in $\calH^{ \times}$ is simply
\begin{equation} 
\frac{1}{a+bj} = \frac{a-bj}{a^2-b^2} 
\end{equation}
this follows from the identity $(a+bj)(a-bj) = a^2-b^2$ given $a^2-b^2 \neq 0$. Let $\Bcal = \RN \times \RN$ with $(a,b)(c,d) = (ac, bd)$ for all $(a,b),(c,d) \in \Bcal$. We can show that
\begin{equation} \Psi(a,b) = a\left(\frac{1+j}{2}\right)+b\left(\frac{1-j}{2}\right) \qquad \& \qquad
\Psi^{-1}(x+jy) = (x+y,x-y) 
\end{equation} 
provide an isomorphism of $\calH$ and $\RN \times \RN$. In \cite{cookAcalculusI} an examples are given which show how this isomorphism can be used to solve the quadratic equation in $\Hcal$ and to derive  d'Alembert's solution to the wave equation.
\end{exa}

\subsection{Submultiplicative norms}
The division algebras $\RN, \CN$ and $\HN$ can be given a {\bf multiplicative norm} where $\| x \star y \| = \| x \| \, \| y \|$. Generally we can only find {\bf submultiplicative norm}. 

\begin{exa}
 If $\Hcal$ is given norm $\| x+jy \| = \sqrt{x^2+y^2}$ then $\| zw \| \leq \sqrt{2} \|z \| \, \|w \|$.
\end{exa} 

\noindent
If $\Acal$ is an algebra over $\RN$ with basis $\{ v_1, \dots , v_n \}$ then define {\bf structure constants} $C_{ijk}$ by $v_i \star v_j = \sum_{k=1}^n C_{ijk} v_k$ for all $1 \leq i,j  \leq n$. For proof of what follows see \cite{cookAcalculusI}.   

\begin{thm} \label{thm:submultiplicative}
\textbf{(submultiplicative norm)} If $\Acal$ is an associative $n$-dimensional algebra over $\RN$ then there exists a norm $|| \cdot ||$ for $\Acal$ and $\Abound >0$ for which $|| x \star y|| \leq \Abound ||x|| ||y||$ for all $x,y \in \Acal$. Moreover, for this norm we find $\Abound = \mathbf{C}(n^2-n+1) \sqrt{n}$ where $\mathbf{C} = \text{max} \{ C_{ijk} \ | \ 1 \leq i,j,k \leq n \}$.
\end{thm}

\begin{coro} \label{prop:ineqnpower}
If $\| x \star y \| \leq \Abound \| x \| \|y \|$ for $x,y \in \Acal$ then $\| z^n \| \leq \Abound^n \| z \|^n$ for each $z \in \Acal$ and $n \in \NN$.
\end{coro}

\begin{coro} \label{thm:quotientinequality}
Suppose $\Abound >0$ is a real constant such that $|| x \star y|| \leq \Abound ||x|| ||y||$ for all $x,y \in \Acal$. If $b \in \Acalx$ and $a \in \Acal$ then $ \frac{||a||}{||b||} \leq \Abound \, \big{|}\big{|} \frac{a}{b} \big{|}\big{|}$. 
\end{coro}

\subsection{Differential calculus on $\Acal$}
The definition of differentiability with respect to an algebra variable is open to some debate. There seem to be two main approaches:
\begin{quote}
\begin{enumerate}[{\bf D1:}]
\item define differentiability in terms of an algebraic condition on the differential,
\item define differentiability in terms of a deleted-difference quotient
\end{enumerate}
\end{quote}
In \cite{cookAcalculusI} it is shown that these definitions are interchangeable on an open set in the context of a commutative semisimple algebra. However, it is also shown that in there exist {\bf D1} differentiable functions which are nowhere {\bf D2}. Hence, we prefer to use {\bf D1} as it is more general. Following \cite{cookAcalculusI} we define differentiability with respect to an algebra variable as follows:

\begin{de} \label{defn:Adiff}
Let $U \subseteq \Acal$ be an open set containing $p$. If $f: U \rightarrow \Acal$ is a function then we say $f$ is {\bf $\Acal$-differentiable at $p$} if there exists a linear function $d_pf \in \EndA$ such that
\begin{equation} \label{eqn:frechetquotAcal}
 \lim_{h \rightarrow 0}\frac{f(p+h)-f(p)-d_pf(h)}{||h||} = 0. 
\end{equation}
 \end{de}

\noindent
In other words, $f$ is $\Acal$-differentiable at a point if its differential at the point is a right-$\Acal$-linear map. Equivalently, given a choice of basis, $f$ is $\Acal$-differentiable if its Jacobian matrix is found in the matrix regular representation of $\Acal$. If $\Acal$ has basis $\beta = \{ v_1, \dots , v_n\}$ has coordinates $x_1, \dots , x_n$ then $d_pf(e_j) = \frac{\partial f}{\partial x_j}(p)$. Suppose $v_1=\mathds{1}$ then $ d_pf(1) = \frac{\partial f}{\partial x_1}(p)$. Observe right linearity of the differential indicates $d_pf(v_j) = d_pf(\mathds{1} \star v_j) = d_pf(1) \star v_j$ hence for each $p$ at which $f$ is $\Acal$ differentiable we find:
\begin{equation}
\frac{\partial f}{\partial x_j} (p)= \frac{\partial f}{\partial x_1}(p) \star v_j.
\end{equation}
These are the {\bf $\Acal$-Cauchy Riemann Equations}. There are $n-1$ equations in $\Acal$ which amount to $n^2-n$ scalar equations. If the $\Acal$-CR equations hold for a continuously differentiable $f$ at $p$ then we have that $d_pf \in \EndA$.  \\

\noindent
Next we wish to explain how to construct the derivative function $f'$ on $\Acal$.  We are free to use the isomorphism between the right $\Acal$ linear maps and $\Acal$ as to define the {\it derivative at a point} for via $f'(p) = \#(d_pf)$. This is special to our context. In the larger study of real differentiable functions on an $n$-dimensional space no such isomorphism exists and it is not possible to identify arbitrary linear maps with points.

\begin{de} \label{defn:derivative}
Let $U\subseteq \Acal$ be an open set and $f: U  \rightarrow \Acal$ an $\Acal$-differentiable function on $U$ then we define $f': U \rightarrow \Acal$ by $f'(p) = \# (d_pf) $ for each $p \in U$.
\end{de}

\noindent
Equivalently, we could write $f'(p) = d_pf( \mathds{1})$ since $\#(T) = T( \mathds{1})$ for each $T \in \EndA$.
Many properties of the usual calculus hold for $\Acal$-differentiable functions.  

\begin{prop}
For $f$ and $g$ both $\Acal$-differentiable at $p$,
 \begin{enumerate}[{\bf (i.)}]
\item $ \ds (f+g)'(p) = f'(p)+g'(p)$,
\item for constant $c \in \Acal$, $ \ds (c \star f)'(p) = c \star f'(p)$,
\item Given $\Acal$ is commutative, $ (f \star g)'(p) = f'(p) \star g(p)+ f(p) \star g'(p)$.
\item $ (f \comp g)'(p) = f'(g(p)) \star g'(p). $
\item If $f(\zeta) = \zeta^n$ for some $n \in \NN$ then $f'(\zeta) = n \zeta^{n-1}$. 
\end{enumerate}
\end{prop}

\noindent
If $\Acal$ is not commutative then the product of $\Acal$-differentiable functions need not be $\Acal$-differentiable. In \cite{cookAcalculusI} an example is given where $f,g$ and $f \star g$ are $\Acal$-differentiable yet $g \star f$ is not $\Acal$-differentiable. \\

\noindent
We are also able to find an $\Acal$-generalization of Wirtinger's calculus. In \cite{cookAcalculusI} we introduce conjugate variables $\bar{\zeta}_2, \dots, \bar{\zeta}_n$ for $\Acal$ and find for commutative algebras if $f: \Acal \rightarrow \Acal$ is $\Acal$-differentiable at $p$ then $\ds \frac{\partial f}{\partial \overline{\zeta}_j} = 0$ for $j=2, \dots , n$. In other words, another way we can look at $\Acal$-differentiable functions is that they are functions of $\zeta$ alone.  \\  

\noindent
The theory of higher derivatives is also developed in \cite{cookAcalculusI}.

\begin{de} \label{defn:higherderivative}
Suppose $f$ is a function on $\Acal$ for which the derivative function $f'$ is $\Acal$-differentiable at $p$ then we define
$ f''(p) = (f')'(p)$. Furthermore, supposing the derivatives exist, we define $f^{(k)}(p) = (f^{(k-1)})'(p)$ for $k =2,3, \dots$. 
\end{de}

\noindent
Naturally we define functions $f'', f''', \dots, f^{(k)}$ in the natural pointwise fashion for as many points as the derivatives exist. Furthermore, with respect to $\beta = \{ v_1, \dots , v_n \}$ where $v_1 = \mathds{1}$, we have $f'(p) = d_pf(\mathds{1}) = \frac{\partial f}{\partial x_1}(p)$.  Thus, $f' = \frac{\partial f}{\partial x_1}$. Suppose $f''(p)$ exists. Note,
\begin{equation}
 f''(p) = (f')'(p) = \#( d_p f'(\mathds{1}) )= \frac{\partial f'}{\partial x_1}(p) = \frac{\partial^2f}{\partial x_1^2}(p). 
\end{equation}
Thus, $f'' = \frac{\partial^2 f}{\partial x_1^2}$. By induction, we find if $f^{(k)}$ exists then $f^{(k)} = \frac{\partial^k f}{\partial x_1^k}$. Furthermore, if $f: \Acal \rightarrow \Acal$ is $k$-times $\Acal$-differentiable then 
\begin{equation}
\frac{\partial^k f}{\partial x_{i_1}\partial x_{i_2} \cdots \partial x_{i_k}} = \frac{\partial^k f}{\partial x_1^k} \star v_{i_{1}} \star v_{i_{2}} \star \cdots \star v_{i_{k}}.
\end{equation}

\noindent
The Theorem below gives us license to convert equations in $\Acal$ to partial differential equations which every component of an $\Acal$-differentiable function must solve!

\begin{thm}
Let $U$ be open in $\Acal$ and suppose $f: U \rightarrow \Acal$ is $k$-times $\Acal$-differentiable. If there exist $B_{i_1i_2\dots i_k} \in \RN$ for which $\sum_{i_1i_2\dots i_k} B_{i_1i_2\dots i_k}v_{i_1} \star v_{i_2} \star \cdots \star v_{i_k} = 0$ then 
\begin{equation}
\sum_{i_1i_2\dots i_k} B_{i_1i_2\dots i_k}\frac{\partial^k f}{\partial x_{i_1} \partial x_{i_2} \cdots  \partial x_{i_k} }=0. 
\end{equation}
\end{thm}

\begin{exa}
Since $i^2=-1$ in $\CN$ it follows for $z=x+iy$ that complex differentiable $f$ have $f_{yy}=-f_{xx}$. We usually see this in notation $f=u+iv$ and the observation $u_{xx}+u_{yy}=0$ and $v_{xx}+v_{yy}=0$. The real and imaginary parts of a complex differentiable function are harmonic because $1+i^2=0$.
\end{exa}

\begin{exa}
Since $j^2=1$ in $\Hcal$ it follows for $z=x+jy$ that $\Hcal$-differentiable $f$ have $f_{yy}=f_{xx}$. If $f=u+iv$ then $u_{xx}-u_{yy}=0$ and $v_{xx}-v_{yy}=0$. Thinking of $y$ as time and $x$ as position, the partial differential equation $u_{xx}=u_{yy}$ is a unit-speed wave equation. Component functions of hyperbolic differentiable functions are solutions to the wave equation in one dimension!\end{exa}

\noindent
For complex algebras the Cauchy Integral Formula links differentiation and integration in such a way that one complex derivative's existence requires all higher complex derivative's likewise exist. We consider many examples where $\Acal$ does not permit such a simplification. However, if we know $f$ is smooth in the real sense and once $\Acal$-differentiable then it is shown in \cite{cookAcalculusI} that $f$ allows infinitely many $\Acal$-derivatives. Taylor's Formula for $\Acal$ is also given:

\begin{thm} \label{thm:AcalTaylor}
\textbf {(Taylor's Formula for $\Acal$-Calculus:)} Let $\Acal$ be a commutative, unital, associative algebra over $\RN$. If $f$ is real analytic at $p \in \Acal$ then
$$  f(p+h) = f(p)+f'(p) \star h + \frac{1}{2}f''(p) \star h^2+ \cdots + \frac{1}{k!}f^{(k)}(p) \star h^k+ \cdots $$
where $h^2 = h \star h$ and $h^{k+1} = h^{k} \star h$ for $k \in \NN$.
\end{thm}

\subsection{Integral calculus on $\Acal$}

\noindent
Integration along curves in $\Acal$ is defined in \cite{cookAcalculusI} in much the same fashion as $\CN$. If $\zeta: [t_o,t_1] \rightarrow \Acal$ is differentiable parametrization of a curve $C$ and $f$ is continuous near $C$ then 
\begin{equation}
 \ds \int_C f ( \zeta) \star d\zeta = \int_{t_o}^{t_f} f( \zeta (t)) \star \frac{d\zeta}{dt} \, dt. \end{equation}

\begin{thm} \label{thm:subML}
Let $C$ be a rectifiable curve with arclength $L$. Suppose $||f(\zeta) || \leq M$ for each $\zeta \in C$ and suppose $f$ is continuous near $C$. Then 
$$ \bigg{|}\bigg{|} \int_C f( \zeta) \star d\zeta \bigg{|}\bigg{|} \leq \Abound ML $$
where $\Abound$ is a constant such that $|| z \star w || \leq \Abound ||z|| \, ||w||$ for all $z,w \in \Acal$.
\end{thm}

\noindent
Let us conclude with a list of notable results given in \cite{cookAcalculusI}. The order in which these results are derived is perhaps surprising. In fact, the next result is last:

\begin{thm} \label{thm:FTCIforalg}
\textbf {(Fundamental Theorem of Calculus Part I:)} Let $C$ be a differentiable curve from $\zeta_o$ to $\zeta$ in $U \subseteq \Acal$ where $U$ is an open simply connected subset of $\Acal$. Assume $f$ is $\Acal$ differentiable on $U$ then 
$$ \frac{d}{d\zeta}\int_C f( \eta) \star d \eta = f(\zeta). $$
\end{thm}

\begin{thm}\textbf {(Fundamental Theorem of Calculus Part II:)} \label{thm:FTCforalg}
Suppose $f = \frac{dF}{d\zeta}$ near a curve $C$ which begins at $P$ and ends at $Q$ then
$$ \int_C f( \zeta) \star d \zeta = F( Q) - F(P). $$
\end{thm}


\begin{thm} \label{thm:topologicalAintegral}
Let $f: U \rightarrow \Acal$ be a function where $U$ is a connected subset of $\Acal$ then the following are equivalent:
\begin{enumerate}[{\bf (i.)}]
\item $\int_{C_1} f \star d\zeta = \int_{C_2} f \star d\zeta$ for all curves $C_1,C_2$ in $U$ beginning and ending at the same points,
\item $\int_C f \star d\zeta =0$ for all loops in $U$,
\item $f$ has an antiderivative $F$ for which $\frac{dF}{d\zeta} = f$ on $U$.
\end{enumerate}
\end{thm}

\begin{thm} \label{coro:adiffloopszero}
\textbf {(Cauchy's Integral Theorem for $\Acal$:)} If $U \subseteq \Acal$ is simply connected then $\int_C f \star d \zeta = 0$ for all loops $C$ in $U$ if and only if $f$ is $\Acal$-differentiable on $U$.
\end{thm}

\section{Sequences} \label{sec:sequences}
In this section we will discuss sequences and the concept of limits and convergence in for a real associative algebra of finite dimension.  To do this, we will generalize many of the theorems from real analysis to our context.  Much of our generalization parallels Rudin's \textit{Principles of Mathematical Analysis} \cite{Rudin}. 
\begin{de}\label{defn:sequence}
A function $f: \NN \rightarrow \Acal$ is called a {\bf sequence} in $\Acal$.  If $f(n) = x_n$, for $n \in \NN$, then it is customary to denote the sequence $f$ by the symbol $\{x_n\}$.  The values of $f$, that is, the elements $x_n$, are called the {\bf terms} of the sequence.   
\end{de}
\noindent
Convergence of sequences is measured in terms of the norm $|| \cdot ||$ on $\Acal$. 
\begin{de}\label{defn:convergence}
Let $\Acal$ be a real associative algebra with a norm $\|\cdot \|$.  A sequence $\{p_n\}$ in $\Acal$ is said to converge if there is a point p $\in\Acal$ with the following property:
For every $\varepsilon>0$ there is an integer M such that $n \geq M$ implies that 
$$\|p_n- p\|<\varepsilon.$$  
In this case we also say that $\{p_n\}$ converges to p, or that p is the limit of $\{p_n\}$, and we write $p_n\to$p, or
$ \ds \lim_{n\to\infty} p_n=p$. If $\{p_n\}$ does not converge, then it is said to diverge.
\end{de}
\begin{exa}
Consider the sequence $\{p_n\}$ in $\CN$ defined by $p_n=\frac{in+1+i}{n}$ for all $n \in \NN.$  Recall $| \cdot |$ defined by $|x+iy| = \sqrt{x^2+y^2}$ provides a norm on $\CN$. Consider,
\begin{equation}
|p_n-i|=\left|\frac{in+1+i}{n}-i\right|=\frac{|in+1+i-in|}{n}=\frac{|1+i|}{n}=\frac{\sqrt{2}}{n}.
\end{equation}
But, given $\varepsilon>0$, we know from the Archimedian property of real numbers that there exists $M \in \NN$ such that $\frac{1}{M}<\frac{\varepsilon}{\sqrt{2}}$.  Thus, for all $n\geq M$, we have:
\begin{equation}
|p_n-i|=\frac{\sqrt{2}}{n}\leq\frac{\sqrt{2}}{M}<\sqrt{2}\left(\frac{\varepsilon}{\sqrt{2}}\right)=\varepsilon.
\end{equation}
Thus $\{ p_n \}$ converges to $i$.
\end{exa}
\noindent
Since $\Acal$ is a vector space, a sequence in $\Acal$ is a sequence of vectors over $\RN$. Part (iii.) of the Theorem below explains that the convergence of a vector sequence is tied to the convergence of its component sequences relative to a basis. In contrast, Parts (i.), (ii.) and (iv.) of the Theorem below directly resemble the usual results for real sequences,the linearity and multiplicativity of limits, following Theorem 3.3 of \cite{Rudin}.

\begin{thm}\label{thm:limit laws}
Suppose $\Acal$ is an associative algebra paired with a submultiplicative norm $\|\cdot\|$ and a basis $\{v_1, \dots ,v_N \}$ such that $\| v_i \|=1$ and for all $x=\sum_{i=1}^{N}x^iv_i, |x^i|\leq\|x\|$ for all $i$.  Suppose also that $\{s_n\}, \{t_n\}$ are sequences in $\Acal$ with $s_n=\sum_{i=1}^{N} s_n^i v_i, \  t_n = \sum_{j=1}^{N} t_n^j v_j$, and $\lim_{n\to\infty} s_n=s, \lim_{n\to\infty} t_n=t$ where $s=\sum_{i=1}^{N}s^i, t=\sum_{j=1}^{N}t^j$.  Then:
\begin{quote}
\begin{enumerate}[{\bf (i.)}]
\item $\lim_{n\to\infty} (s_n+t_n)=s+t$,
\item $ \lim_{n\to\infty} (\alpha\star s_n)=\alpha\star s$, for any number $\alpha\in\Acal$,
\item $ \lim_{n\to\infty}s_n=s$  if and only if $lim_{n\to\infty}s_n^i=s^i$ for all 
$ i =1,2,\dots , N$, 
\item $ \lim_{n\to\infty} (s_n\star t_n)=s\star t$.
\end{enumerate}
\end{quote}
\end{thm}

\noindent
\textbf {Proof:} begin with {\bf (i.)} and {\bf (ii.)} suppose $s_n \rightarrow s$ and $t_n \rightarrow t$ with respect to $|| \cdot ||$ on $\Acal$ as described in the Theorem. For $\alpha\neq 0$, given $\varepsilon>0$, there exist integers $N_1, N_2$ such that
\begin{align} 
&n\geq N_1\implies \|s_n-s\|<\frac{\varepsilon}{2\Abound\|\alpha\|},\\ \notag
&n\geq N_2\implies \|t_n-t\|<\frac{\varepsilon}{2}.
\end{align}
If $N_3= \text{max} \{ N_1, N_2 \}$ then $n \geq N_3$ implies
\begin{align}
\|(\alpha\star s_n+t_n)-(\alpha\star s+t)\|&\leq\|\alpha\star(s_n-s)\|+\|t_n-t\|\\ \notag
&\leq \Abound\|\alpha\|\|(s_n-s)\|+\|t_n-t\|\\ \notag
&<\Abound\|\alpha\|\frac{\varepsilon}{2\Abound\|\alpha\|}+\frac{\varepsilon}{2}\\ \notag
&=\varepsilon.
\end{align}
Let $\alpha=1$ to obtain (i.) and let $t_n=0$ for all $n\in\NN$ to obtain (ii.) for $\alpha \neq 0$. If $\alpha=0$, then for any $\varepsilon >0$ we note $|| \alpha\star s_n - \alpha \star s|| = ||0-0||= 0$ hence $ \lim_{n\to\infty} (\alpha\star s_n)=\alpha\star s$. \\

\noindent
Next we give the proof of {\bf(iii.):} \\
($\Rightarrow$) Assume $s_n \rightarrow s$. For $\varepsilon>0$, suppose there exists $M\in\NN$ such that $n\geq M$ implies $\|s_n-s\|<\varepsilon$. Thus $n\geq M$ implies 
$ |s_n^i-s^i|\leq \| s_n-s \|<\varepsilon$ and we find $s_n^i \rightarrow s^i$ for all $i=1, \dots N$. \\

\noindent
($\Leftarrow$) Assume $s_n^i \rightarrow s^i$ for all $i=1,2,\dots , N$. Given $\varepsilon>0$, choose $ M_1, M_2,...M_N \in \NN$ such that $n\geq M_i$ implies $|s_n^i-s^i|<\frac{\varepsilon}{N}$.
If $M=\text{max}\{M_1, \dots , M_N \}$ and $n \geq M$ then we find
\begin{align}
\|s_n-s\|&=\left\|\sum_{i=1}^{N}(s_n^i-s^i)v_i \right\|\\ \notag
&\leq\sum_{i=1}^{N}\|(s_n^i-s^i)v_i\|\\ \notag
&=\sum_{i=1}^{N}|s_n^i-s^i|\| v_i\|\\ \notag
&<\sum_{i=1}^{N} \frac{\varepsilon}{N}(1)  \\ \notag
&=\varepsilon.
\end{align}
Therefore, $s_n \rightarrow s$ and this completes the proof of {\bf (iii.)}. \\

\noindent
The proof {\bf (iv.)} requires some calculation. Let $C_{ijk}$ be constants such that $v_i \star v_j = \sum_k C_{ijk} v_k$ and $\mathbf{C} = \text{max} \{ |C_{ijk}| \ | \ 1 \leq i,j,k \leq n \}$. Applying the triangle inequality we find:
\begin{equation} \label{eqn:ncestimate}
\| v_i \star v_j \| \leq \sum_{k=1}^N \| C_{ijk}v_k \|  \leq \sum_{k=1}^N   \mathbf{C} \| v_k \|  =\mathbf{C} \| v_k \|\sum_{k=1}^N 1 =  \mathbf{C}N
\end{equation}
as $\| v_k \| = 1$ and $\sum_{k=1}^N 1 = N$. Calculate:
\begin{align} \label{eqn:stackoestimatesI}
\|s_n\star t_n-s\star t\|
&=\bigg{\|}\sum_{i,j}(s_n^i v_i)\star (t_n^j v_j)-\sum_{i,j}(s^i v_i)\star (t^j v_j) \bigg{\|}\\ \notag
&= \bigg{\|}\sum_{i,j}(s_n^it_n^j-s^it^j)(v_i\star v_j) \bigg{\|}\\ \notag
&\leq \sum_{i,j}\left|s_n^it_n^j-s^it^j\right| \| v_i\star v_j \|  \\ \notag
&\leq \sum_{i,j}\left|s_n^it_n^j-s^it^j\right|\mathbf{C}N 
\end{align}
where we have used Equation \ref{eqn:ncestimate} in the last step. 
Observe that
\begin{equation}
s_n^it_n^j-s^it^j=(s_n^i-s^i)(t_n^j-t^j)+s^i(t_n^j-t^j)+t^j(s_n^i-s^i)
\end{equation}
and apply it to Equation \ref{eqn:stackoestimatesI} as to obtain:
\begin{align} \label{eqn:stackoestimates}
\|s_n\star t_n-s\star t\|
&= \mathbf{C}N\sum_{i,j} \left|(s_n^i-s^i)(t_n^j-t^j)+s^i(t_n^j-t^j)+t^j(s_n^i-s^i)\right|\\ \notag
&\leq \mathbf{C}N\biggl(\sum_{i,j}|s_n^i-s^i||t_n^j-t^j|+\sum_{i,j}|s^i||t_n^j-t^j|+\sum_{i,j}|t^j||s_n^i-s^i|\biggr)\\ \notag
&\leq \mathbf{C}N\biggl(\sum_{i,j}|s_n^i-s^i||t_n^j-t^j|+\sum_{i,j}M_s|t_n^j-t^j|+\sum_{i,j}M_t|s_n^i-s^i|\biggr) \\ \notag
&\leq \mathbf{C}N\biggl(\sum_{i,j}|s_n^i-s^i||t_n^j-t^j|+\sum_{j}NM_s|t_n^j-t^j|+\sum_{i}NM_t|s_n^i-s^i|\biggr)
\end{align}
where $M_s=\text{max}\{s^1, \dots, s^N\}$ and  $M_t=\text{max}\{t^1, \dots, t^N\}$.  Let $\varepsilon>0$. Assume $s_n \rightarrow s$ and $t_n \rightarrow t$ thus by ({\bf iii.}) we know the component sequences $s_n^i \rightarrow s^i$ and $t_n^i \rightarrow t^i$. It follows we may choose $N_1^i, N_2^i, N_3^i, N_4^i \in \NN$ for $i = 1, \dots , N$ such that
 \begin{align} \label{eqn:estimotto}
 &n\geq N_1^i \ \ \Rightarrow \ \ |s_n^i-s^i|<\frac{\sqrt{\varepsilon N}}{N^2\sqrt{2\mathbf{C}}}, \\ \notag
 &n\geq N_2^i\ \ \ \Rightarrow \ \ |t_n^i-t^i|<\frac{\sqrt{\varepsilon N}}{N^2\sqrt{2\mathbf{C}}},\\ \notag
 &n\geq N_3^i\ \ \ \Rightarrow \ \ |s_n^i-s^i|<\frac{\varepsilon}{4N^2M_s\mathbf{C}}, \\ \notag
 &n\geq N_4^i\ \ \ \Rightarrow \ \ |t_n^i-t^i|<\frac{\varepsilon}{4N^2M_t\mathbf{C}}.
 \end{align}
 If $N_5=\text{max}\{ N_j^i \ | \ j=1,2,3,4, i = 1, \dots , N  \}$ then following Equation \ref{eqn:stackoestimates} and \ref{eqn:estimotto} we find
\begin{align} 
 \|s_n\star t_n-s\star t\|&\leq \mathbf{C}N\biggl( \sum_{i,j=1}^{N}\frac{\sqrt{\varepsilon N}}{N^2\sqrt{2\mathbf{C}}} \cdot \frac{\sqrt{\varepsilon N}}{N^2\sqrt{2\mathbf{C}}}+\sum_{j=1}^{N}\frac{NM_s\varepsilon}{4N^2M_s\mathbf{C}}+\sum_{i=1}^{N}\frac{NM_t\varepsilon}{4N^2M_t\mathbf{C}}\biggr) \\ \notag
&=  \frac{\varepsilon}{2N^2}\sum_{i,j=1}^{N}1+\frac{\varepsilon}{4N}\sum_{j=1}^{N}1+\frac{\varepsilon}{4N}\sum_{i=1}^{N}1 \\ \notag
 &=\varepsilon.
\end{align}
Therefore, $s_n \star t_n \rightarrow s \star t$. $\square$\\

\noindent
Continuity of functions on $\Acal$ can be described sequentially.

\begin{thm}\label{thm:continuity}
The following are equivalent:
\begin{enumerate}[{\bf (i.)}]
\item For all $\varepsilon>0$, there exists $ \delta>0$ such that $0<\|z-z_0\|<\delta$ implies $\|f(z)-f(z_0)\|<\varepsilon$.
\item For all $\{z_n\}$ such that $z_n\rightarrow z_0\ as\  n\rightarrow\infty,\ \lim_{n \to \infty}f(z_n)=f(z_0)$.
\end{enumerate}
\end{thm}

\noindent
{\bf Proof:} {\bf (i.)} $\Rightarrow$ {\bf (ii.)}: Suppose {\bf (i.)}  and let $\{z_n\}$ be a sequence in $\Acal$ such that $\lim_{n\to\infty}z_n=z_0$. Let $\varepsilon>0$.  Choose $\delta>0$ such that $\|z-z_0\|<\delta$ implies $\|f(z)-f(z_0)\|<\varepsilon$.
As $\{z_n\}$ converges to $z_0$, we know there exists $ M \in \NN$ such that $n>M\ \Rightarrow \|z_n-z_0\|<\delta$.  Observe, if $n>M$ then $\|z_n-z_0\|<\delta$ and thus $\|f(z_n)-f(z)\|<\varepsilon$.  Consequently, $\lim_{n\to\infty}f(z_n)=f(z_0)$.\\

\noindent
{\bf (ii.)} $\Rightarrow$ {\bf (i.)}:  Suppose {\bf (ii.)} and assume towards a contradiction that there exists $\varepsilon_0>0$ such that for any $\delta>0$, there exists $z^*$ such that $\|z^*-z_0\|<\delta$ and $\|f(z^*)-f(z_0)\|\geq\varepsilon_0$.
In particular, this holds for $\delta=1/n$.  For each $n\in\NN$ we choose an element $z_n$ such that $\|z_n-z_0\|<1/n$ and $\|f(z_n)-f(z_0)\|\geq\varepsilon_0$ and thus construct the sequence $\{z_n\}$.  Therefore, $z_n\rightarrow z_0$ and $\lim_{n \to \infty}f(z_n)\neq f(z_0)$. But this is a contradiction to {\bf (ii.)} and conclude that {\bf (i.)} is true.  $\square$

\begin{exa}
On any algebra $\Acal$ with a submultiplicative norm described above, the function $f:\ \Acal\to\Acal$ defined by $f(z)=z^2 = z \star z$ for all $z\in\Acal$ is continuous.  Indeed, for any $z_0\in\Acal$ and any sequence $\{z_n\}$ which converges to $z_0$, we have: using (iv.) of Theorem \ref{thm:limit laws},
\begin{equation}
\lim_{n\to\infty}f(z_n)=\lim_{n\to\infty}z_n \star z_n=(\lim_{n\to\infty}z_n) \star (\lim_{n\to\infty}z_n)=z_0\star z_0=z_0^2=f(z_0)\end{equation}
which, by Theorem 3.5, shows continuity at every point of $\Acal$ and thus continuity on $\Acal$.
\end{exa}

\noindent
We close this section with a valuable concept, based on Definition 3.8 of \cite{Rudin}.

\begin{de}\label{def:Cauchy sequence}
A sequence $\{p_n\}$ in an algebra $\Acal$ with a norm $\|\cdot\|$ is said to be a {\bf Cauchy sequence} if for every $\varepsilon>0$ there is an integer $M$ such that if $n,\ m\geq M$ then $\|p_n-p_m\|<\varepsilon$.
\end{de}
\begin{thm}\label{thm:Complete}
\textbf{(Cauchy Criterion)}  If $\Acal$ is an associative algebra paired with a submultiplicative norm $\|\cdot\|$ and a basis $\{ v_1, \dots , v_N\}$ such that for all $x=\sum_{i=1}^{N}x^i v_i, |x^i|<\|x\| $ for each $i=1, \dots, N$, then a sequence $\{p_n\}$ in $\Acal$ is is convergent if and only if it is Cauchy.
\end{thm}

\noindent
\textbf {Proof:} ($\Rightarrow$) Suppose the sequence $\{p_n\}$ converges to p.  Given $\varepsilon>0$,  there exists an integer $ M>0$ such that $n\geq M$ implies $\|p_n-p\|<\frac{\varepsilon}{2}$.  Suppose $m,\ n\geq N$
\begin{equation}
\|p_n-p_m\|=\|p_n-p+(p-p_m)\|\leq\|p_n-p\|+\|p_m-p\|<\frac{\varepsilon}{2}+\frac{\varepsilon}{2}=\varepsilon.
\end{equation}
Thus $\{p_n\}$ is a Cauchy sequence. \\

\noindent
($\Leftarrow$) Suppose the sequence $\{p_n\}$ is Cauchy.  Thus, given $\varepsilon>0$, there exists and integer $M>0$ such that $n,\ m>N$ implies $\|p_n-p_m\|<\varepsilon$. Consider, $p_n-p_m=\sum_{i=1}^{n}(p_n^i-p_m^i)v_i$.  Thus, given  the properties of the basis of $\Acal$, 
$|p_n^i-p_m^i|<\|p_n-p_m\|<\varepsilon$ for $i=1, \dots , N$. We have shown all the component sequences $\{p_n^i\}$ are Cauchy. We know from real analysis that a real sequence converges if and only if it is Cauchy.  Thus, for each $i=1, \dots , N$,  $\{p_n^i\}$ is a real Cauchy sequence and there exists $p^i\in\RN$ for which $p_n^i \rightarrow p^i$. Thus, defining $p = \sum_{i=1}^N p^i v_i$ we find $p_n \rightarrow p$ by part (iii.) of Theorem \ref{thm:limit laws}.  $\square$\\ 

\noindent
If $\Acal$ meets the Cauchy Criterion, we say that it is a complete algebra.
\begin{exa}
$\RN^{n\times n}$ with $\| A \| = \sqrt{ \text{trace}( A^TA) }$ is a complete algebra since its norm has the necessary property with respect to the basis of unit-matrices $\{ E_{ij} \}$ defined by $(E_{ij})_{kl} = \delta_{ik}\delta_{jl}$. It is not hard to show $|A_{ij}| \leq \| A \|$ and $\| E_{ij} \| =1$ for all $i,j$.
\end{exa}
\noindent
All the examples we study are complete since it is known that any finite dimensional vector space over $\RN$ is complete. For example, see Theorem 2.4-2 on page 73 of \cite{Kreyszig}.

\section{Series} \label{sec:series}
We now consider sequences of sums and the limits of these sequences.  The definition and the two theorems that follow mirror Rudin's 3.21-3.23 in \cite{Rudin}
\begin{de}\label{def:series}
Given a sequence $\{a_n\}$ in an algebra $\Acal$, we call the sequence $\{s_n\}$ where $ \ds s_n=\sum_{k=1}^{n}a_k$ the {\bf sequence of partial sums} of $\{a_n\}$.  We denote $\ds \lim_{n\to\infty}s_n=\sum_{k=1}^{\infty}a_k$ or  $\sum a_k$ and call this limit an infinite series. If $\{s_n\}$ converges to s, we say that the {\bf series converges} and write $\sum_{k=1}^{\infty}a_k=s$. If $\{s_n\}$ diverges, we say that the {\bf series diverges}.
\end{de}

\noindent
Notice if $m \geq n$ and $s_n=\sum_{k=1}^{n}a_k$ then
\begin{equation}
 s_m-s_n = \sum_{k=1}^{m}a_k - \sum_{k=1}^{n}a_k = \sum_{k=m}^{n}a_k. 
\end{equation}
Hence the Cauchy criterion (Theorem \ref{thm:Complete}) applied to the sequence of partial sums provides the following useful test for convergence of a series:
\begin{thm}\label{thm:Cauchy series}
(\textbf{ Cauchy criterion for series}) The series $\sum a_k$ in $\Acal$ converges if and only if for every $\varepsilon>0$ there exists $M \in \NN$ such that $m\geq n\geq M$ implies $ \ds \bigg{\|} \sum_{k=n}^{m}a_k \bigg{\|} \leq\varepsilon$.
\end{thm}

\noindent
In particular, if $\sum a_k$ converges then by taking $m=n$, we find $\|a_n\|\leq\varepsilon$ for all $n \geq M$. Therefore, we find:

\begin{thm}\label{thm:term test}
\textbf{($n$-th Term Test)} If $\sum a_n$ converges then $\lim_{n\to\infty} a_n=0$.
\end{thm}

\noindent
We now develop some theorems to test for convergence of series.  The next theorem follows 3.25 of \cite{Rudin}, with a modified second part.
\begin{thm}\label{thm:comparison}
\textbf{(Comparison Test)}  If $\|a_n\|\leq c_n$ for $ n\geq N_0$, where $N_0$ is some fixed integer, and if $\sum c_n$ converges, then $\sum a_n$ converges. Likewise, if there exists $b_n \in [0, \infty)$ for which $\sum_n b_n$ diverges and $b_n \leq \| a_n \|$ then $\sum \| a_n \|$ diverges.
\end{thm}

\noindent
{\bf Proof:} Suppose $\varepsilon>0$. There exists $M\geq N_0$ such that $m\geq n\geq M$ implies $\sum_{k=n}^{m}c_k\leq\varepsilon$, by the Cauchy criterion.  Hence
\begin{equation}
\bigg{\|} \sum_{k=n}^{m}a_k\bigg{\|} \ \leq\ \sum_{k=n}^{m}\|a_k\|\ \leq\ \sum_{k=n}^{m}c_k\ \leq\ \varepsilon. 
\end{equation}
Therefore $\sum a_k$ converges as we have shown it satisfied  Cauchy criterion for series. For the divergent case, since $\sum b_n$ diverges it follows its sequence of partial sums are unbounded and hence $\| a_n \|$ also has a unbounded sequence of partial sums hence $\sum_n \| a_n \|$ diverges. $\square$
\begin{exa}
In $\RN^{2 \times 2}$ where $\left \| \left(\begin{array}{cc} a & b \\ c & d \end{array}\right)\right\| = \sqrt{a^2+b^2+c^2+d^2}$, the series $\ds \sum_{n=1}^\infty \left[ \begin{array}{rr} 1/n^2 & 2/n^2 \\ 3/n^2 & 4/n^2 \end{array}\right]$ converges. For every $n$, we have $\left\|\left( \begin{array}{rr} 1/n^2 & 2/n^2 \\ 3/n^2 & 4/n^2 \end{array}\right)\right\|=\frac{\sqrt{1+4+9+16}}{n^2}=\frac{\sqrt{30}}{n^2}$.  But, we know from real analysis that $\sum_{n=1}^\infty\frac{
\sqrt{30}}{n^2}$ converges, so the result follows by the Comparison Test. 
\end{exa}

\begin{de}
If the series $\sum \|a_n\|$ converges then $\sum a_n$ is said to {\bf converge absolutely}.
\end{de}

\noindent
Naturally, we recover the standard meaning of absolute convergence for real series given the choice $\Acal = \RN$ with $ \| x \| = \sqrt{x^2}$. 

\begin{exa}
In $\CN$, the series $\sum\frac{i^n}{n}$ converges, as can be shown in complex analysis, but it does not converge absolutely.  Indeed, $\sum|\frac{i^n}{n}|=\sum\frac{1}{n}$ does not converge.
\end{exa}

\begin{thm}
If $\sum \| a_n \|$ converges then $\sum a_n$ converges. In other words, absolute convergence implies ordinary convergence for series in $\Acal$.
\end{thm}

\noindent
{\bf Proof:} Suppose $\sum \| a_n \|$ converges. Let $n>m$
\begin{equation}
 0 \leq \| \sum^n a_k - \sum^m a_k \| = \| \sum_{k=m}^n a_k \| \leq \sum_{k=m}^n \| a_k \|
\end{equation}
Since $\sum \| a_n \|$ converges we know by Cauchy Criterion $\sum_{k=m}^n \| a_k \| \rightarrow 0$. Hence, 
\begin{equation}
 \| \sum^n a_k - \sum^m a_k \| \rightarrow 0 \end{equation} 
and by Theorem \ref{thm:Cauchy series} we find $\sum a_k$ converges. $\Box$ \\

\noindent
The next two theorems mirror 3.33 and 3.34 of \cite{Rudin}.

\begin{thm}\label{thm:numericalRoottest}
\textbf {(Root Test)} Given $\sum a_n$, put $\alpha=\limsup_{n\to\infty}\sqrt[n]{\|a_n\|}$.  Then:
\begin{quote}
\begin{enumerate}[{\bf (i.)}]
\item if $\alpha<1$, $\sum a_n$ converges absolutely;
\item if $\alpha>1$, $\sum a_n$ diverges;
\item if $\alpha=1$, the test gives no information.
\end{enumerate}
\end{quote}
\end{thm}

\noindent
\textbf {Proof:}  If $\alpha<1$, we can choose $\beta$ such that $\alpha<\beta<1$, and an integer $M$ such that 
\begin{equation}
\sqrt[n]{\|a_n\|}<\beta
\end{equation}
for $n\geq M$.  That is, $n\geq M$ implies
\begin{equation}
\|a_n\|<\beta^n.
\end{equation}
Since $0<\beta<1$ we recognize the geometric series $\sum\beta^n$ converges.   Convergence of $\sum \| a_n \|$ follows from the comparison test (Theorem \ref{thm:comparison}). To prove {\bf (ii.)} supppose $\alpha >1$ then by definition of limsup there exists a subsequence $\sqrt[n_k]{\|a_{n_k}\|} \rightarrow \alpha >1$ thus $a_n \nrightarrow 0$ hence $\sum a_n$ diverges. {\bf (iii.)} the usual examples from real calculus suffice. $\Box$ \\

\begin{exa}
In $\CN$, the series $\sum(\frac{i \cos(n)}{2})^n$ converges by the Root Test.  Indeed, we have:
\begin{equation}
\alpha=\limsup_{n\to\infty}\sqrt[n]{\left|\frac{i\cos(n)}{2}\right|^n}=\limsup_{n\to\infty}\left|\frac{i \cos(n)}{2}\right|=\limsup_{n\to\infty}\frac{\cos(n)}{2}=\frac{1}{2}<1. 
\end{equation}
\end{exa}

\begin{thm}\label{thm:Ratio test}
\textbf{(Ratio Test)}  The series $\sum a_n$
\begin{quote}
\begin{enumerate}[{\bf (i.)}]
\item converges absolutely if $\limsup_{n\to\infty}\frac{\|a_{n+1}\|}{\|a_n\|}<1$,\par
\item diverges if $\frac{\|a_{n+1}\|}{\|a_n\|}\geq1$ for all $n\geq n_0$, where $n_0$ is some fixed integer.
\end{enumerate}
\end{quote}
\end{thm}

\noindent
{\bf Proof:} If condition {\bf(i.)} holds, we can find $\beta<1$, and an integer $M$, such that
\begin{equation}
\frac{\|a_{n+1}\|}{\|a_n\|}<\beta
\end{equation}
for $n\geq M$.  In particular,
\begin{equation}
\|a_{M+1}\| <\beta \|a_M\| \  \Rightarrow \  
\|a_{M+2}\|<\beta \|a_{M+1}\|<\beta^2 \|a_M\| \  \Rightarrow  \  \cdots \  \Rightarrow \  \|a_{M+p}\|<\beta^p \|a_M\|.
\end{equation}
Thus, $\|a_n\|<\|a_M\|\beta^{-M}\cdot \beta^n$
for $n\geq M$, and it follows from the comparison test that $\sum \|a_n\|$ converges, since $\sum \beta^n$ converges, and thus we obtain {\bf (i.)}. \\

\noindent
To understand {\bf (ii.)} suppose $\|a_{n+1}\|\geq \|a_n\| \neq 0$ for $n\geq n_0$, it is easily seen that the condition $a_n\rightarrow 0$ fails and thus {\bf (ii.)} follows by the $n$-th term test. $\Box$ \\

\noindent
Care should be taken with part (ii.), the absence of a limiting process is significant. The knowledge that lim $a_{n+1}/a_n=1$ implies nothing about the convergence of $\sum a_n$.  The series $\sum 1/n$ and $\sum 1/n^2$ demonstrate this.

\begin{exa}
In the quaternions, the series $\sum\frac{(1+i+j+k)^n}{n!}$ converges absolutely by the Ratio Test.  Note,
\begin{equation}
\alpha=\limsup_{n\to\infty}\frac{\left\|\frac{(1+i+j+k)^{n+1}}{(n+1)!}\right\|}{\left\|\frac{(1+i+j+k)^n}{n!}\right\|}=\limsup_{n\to\infty}\frac{\frac{2^{n+1}}{(n+1)!}}{\frac{2^n}{n!}}=\limsup_{n\to\infty}\frac{2}{n+1} = 0 < 1.
\end{equation}
\end{exa}

\noindent
Finally, we examine the convergence of sums and products of series.
\begin{thm}\label{thm:Sum of series}
If $\sum a_n$ and $\sum b_n$ converge and $c \in \Acal$ then 
$$ \sum (a_n+b_n) = \sum a_n+\sum b_n  \qquad \& \qquad \sum c \star a_n = c \star \sum a_n. $$
\end{thm}

\noindent
{\bf Proof:} Suppose there exist $A,B \in \Acal$ for which $\sum a_n = A$ and $\sum b_n = B$. Partial sums
\begin{equation}
A_n=\sum_{k=0}^{n}a_k \qquad \& \qquad B_n=\sum_{k=0}^{n}b_k
\end{equation}
have $A_n \rightarrow A$ and $B_n \rightarrow B$.
By Theorem \ref{thm:limit laws} we calculate for $c \in \Acal$,
\begin{equation}
c\star A_n+B_n \rightarrow c \star A+B. 
\end{equation}
This proves the theorem. $\Box$ \\

\noindent
The multiplication of two series is understood in terms of the Cauchy product. 

\begin{de}
Given $\sum a_n=A$, and $\sum b_n=B$, we set
$$c_n=\sum_{k=0}^{n}a_k \star b_{n-k} $$
for $n\in\NN$ and call $\sum c_n$ the {\bf product} of the two given series.
\end{de}

\noindent
Absolute convergence is useful to study the existence of the product of two series. In particular, we show that the product of a convergent series with an absolutely convergent series converges to the Cauchy product:

\begin{thm} \label{thm:productofseries}
Suppose
\begin{quote}
\begin{enumerate}[{\bf (i.)}]
\item $\sum_{n=0}^{\infty}a_n\ converges\ absolutely,$
\item  $ \sum_{n=0}^{\infty}a_n=A \in \Acal$ and $\sum_{n=0}^{\infty}b_n=B \in \Acal$,
\item  $ c_n=\sum_{k=0}^{n}a_k \star b_{n-k}$  for all $n\in\NN$.
\end{enumerate}
\end{quote}
Then $\sum_{n=0}^{\infty}c_n=A \star B$.
\end{thm}

\noindent
{\bf Proof:} assume {\bf (i.)}, {\bf (ii.)} and {\bf (iii.)} from the statement of the Theorem are given. Let
\begin{equation}
A_n=\sum_{k=0}^{n}a_k,\ \ \ \ B_n=\sum_{k=0}^{n}b_k,\ \ \ \ C_n=\sum_{k=0}^{n}c_k,\ \ \ \ \beta_n=B_n-B.
\end{equation}
Then notice we can rearrange the terms in the finite sum $C_n$ as follows:
\begin{align}
C_n&=a_0 \star b_0+(a_0 \star b_1+a_1 \star b_0)+\dots+(a_0 \star b_n+a_1 \star b_{n-1}+\dots +a_n \star b_0)\\ \notag
&=a_0 \star B_n+a_1 \star B_{n-1}+\dots+a_n \star B_0\\ \notag
&=a_0 \star (B+\beta_n)+a_1 \star (B+\beta_{n-1})+\dots+a_n \star \beta_0\\ \notag
&=A_n \star B+\underbrace{a_0 \star \beta_n+a_1 \star \beta_{n-1}+\dots+a_n \star \beta_0}_{  \gamma_n}.
\end{align}
We wish to show that $C_n\rightarrow A\star B$.  Since $A_n \star B\rightarrow A\star B$, it suffices to show that $\gamma_n \rightarrow 0$.
Since $\sum a_n$ converges absolutely, we know there exists $\alpha \in [0,\infty)$ for which 
\begin{equation}
\alpha=\sum_{n=0}^{\infty}\|a_n\|.
\end{equation}
Let $\varepsilon>0$ be given.  By {\bf (ii.)}, $\beta_n\rightarrow 0$.  Hence we can choose $M \in \NN$ such that $\|\beta_n\|\leq\frac{\varepsilon}{\alpha\Abound}$ for $n\geq M$, in which case
\begin{align}
\|\gamma_n\|&\leq\|a_n \star \beta_0+\dots+a_{n-M} \star \beta_M\|+\|a_{n-M-1}\star \beta_{M+1}+\dots+a_0 \star \beta_n\|\\ \notag
&\leq\|a_n\star\beta_0+\dots+a_{n-M}\star \beta_M\|+\Abound\|a_{n-M-1} \|\| \beta_{M+1} \|+\dots+\Abound\| a_0 \| \| \beta_n\|\\ \notag
&\leq\|a_n \star \beta_0+\dots+a_{n-M}\star \beta_M\|+\Abound\|a_{n-M-1} \|\frac{\varepsilon}{\alpha\Abound}+\dots+\Abound\| a_0 \| \frac{\varepsilon}{\alpha\Abound} \\ \notag
&\leq\|a_n \star\beta_0+\dots+a_{n-M}\star\beta_M\|+\left( \|a_{n-M-1} \| +\dots+\| a_0 \| \right) \frac{\varepsilon}{\alpha} \\ \notag
&\leq\|a_n\beta_0+\dots+a_{n-M}\beta_M\|+\varepsilon.
\end{align}
Fix $M$ and let $n\rightarrow\infty$, we find
\begin{equation}
\limsup_{n\to\infty}\|\gamma_n\|\leq \varepsilon
\end{equation}
since $a_k\rightarrow 0$ as $k\rightarrow\infty$.  Since $\varepsilon$ is arbitrary, we find $\gamma_n \rightarrow 0$ and the Theorem follows. $\square$\\ 

\noindent
If $\Acal$ is not commutative then it is possible that $A \star B \neq B \star A$. However, it is clear that a similar argument could be given if we were instead given the absolute convergence of $\sum b_n$. Consequently, the product of two convergent series converges to the product of their sums if at least one of the two series converges absolutely.

\section{Power series} \label{sec:powerseries}

\begin{de}
Suppose there exists $z_0 \in \Acal$ and $c_0,c_1, c_2, \dots  \in \Acal$ and
$$ f(z) = \sum_{n=0}^{\infty}c_n \star (z-z_0)^n$$
for each $z \in \Acal$ for which the series converges. Then we say $f(z)$ is a {\bf power series centered at $z_0$ in $\Acal$} with {\bf coefficients} $c_n$. 
\end{de}

\noindent
The domain of a power series is controlled by both its coefficients and its center. If we change the center while holding the coefficients fixed then the domain is modified by translation.  

\begin{lem} \label{lem:centersub}
If $f(z) = \sum_{n=0}^{\infty} c_n \star (z-z_1)^n$ converges on $U \subseteq \Acal$ then $g(z) = \sum_{n=0}^{\infty} c_n \star (z-z_2)^n$ converges on $z_2-z_1+U = \{ z_2-z_1+u \ | \ u \in U \}$.  If $\sum_{n=0}^{\infty} c_n \star z^n$ converges on $U$ then $\sum_{n=0}^{\infty} c_n \star (z-z_0)^n$ converges on $z_0+U$.
 \end{lem}
 
 \noindent
 {\bf Proof:} Suppose $x = z_2-z_1+z$ for $z \in U$ then $x-z_2 = z-z_1$ hence
 \begin{equation} g(x) = \sum_{n=0}^{\infty} c_n \star (z-z_1)^n = f(z). 
 \end{equation}
Since and $f(z)$ exists for $z \in U$ we find $g(x)$ exists for each $x \in z_2-z_1+U$. Finally, set $z_1=0$ and $z_2=z_0$ to obtain the last claim of the Lemma. $\Box$ \\

\noindent
The result below differs from the usual Root Test of real or complex series in that the test does not guarantee divergence for $ \| z \| > R$. 
\begin{thm} \textbf{(Root Test for Power Series)}   \label{thm:roottestpowerseries}
Given an algebra $\Acal$ with $ \| x \star y \| \leq \Abound \|x \| \|y \|$ for all $x,y \in \Acal$ and power series $\sum c_n \star (z-z_o)^n$ in $\Acal$, let
$$\alpha = \limsup_{n\to\infty}\sqrt[n]{\|c_n\|} \qquad \& \qquad  \ R=\frac{1}{ \Abound \alpha }. $$
Then $\sum c_n \star (z-z_o)^n$ is absolutely convergent for $\|z-z_o\|<R$. Moreover, if $\alpha=0$ then $\sum c_n \star (z-z_o)^n$ converges absolutely on $\Acal$.
\end{thm}

\noindent
{\bf Proof:} Let us study power series $\sum c_n \star z^n$ in $\Acal$ centered at $0$.
Let $a_n=c_n \star  z^n$ and seek to apply Theorem \ref{thm:numericalRoottest}.  Consider:
\begin{align}
\limsup_{n\to\infty}\sqrt[n]{\|a_n\|}&=\limsup_{n\to\infty}\sqrt[n]{\|c_n \star z^n \|} \\ \notag
&\leq\limsup_{n\to\infty}\sqrt[n]{m_\Acal\|z^n\| \| c_n\|}\\ \notag
&\leq\limsup_{n\to\infty}\sqrt[n]{m_\Acal^n\|z\|^n\| c_n\|}  \ \ \ \text{(applied Proposition \ref{prop:ineqnpower})}\\ \notag
&=m_\Acal\|z\|\limsup_{n\to\infty}\sqrt[n]{\|c_n\|} \\ \notag
&=\frac{\|z\|}{R}.
\end{align}
Hence, the series is absolutely convergent if $\| z \| < R$. If $\alpha=0$ then Theorem \ref{thm:numericalRoottest} provides absolute convergence at each $z \in \Acal$.  Finally, Lemma \ref{lem:centersub} completes the proof for $z_o \neq 0$.  $\Box$ \\

\noindent
In the theory of power series over the real or complex numbers the root test provides a boundary between points of convergence and divergence. However, we will see in Example \ref{exa:bandconverge} the appearance of zero divisors  makes it is possible to find additional points of convergence beyond those indicated by the root test. 


\begin{thm} \label{thm:norootcriterianeeded}
Suppose that $\sum c_n \star (z-z_o)^n$ converges for all $\|z-z_o\|<R$ for some $R>0$.  Then $\alpha=\limsup_{n\to\infty}\sqrt[n]{\|c_n\|} \leq \frac{\| 1 \|}{R}$.
\end{thm}

\noindent
{\bf Proof:} Consider $z = z_o+\frac{\mathds{1}}{\| \mathds{1} \|}b$ where $0< b < R$ and $\mathds{1}$ denotes the unity in the algebra. Notice $\| z-z_o \| = b  < R$ . Note, if $\alpha=\limsup_{n\to\infty}\sqrt[n]{\|c_n\|}$ then 
\begin{equation}
\limsup_{n\to\infty}\sqrt[n]{\|c_n\star (z-z_o)^n\|} 
 = \frac{b}{\| \mathds{1}\|}\limsup_{n\to\infty}\sqrt[n]{\left\|c_n \right\|} = \frac{b \alpha}{\| \mathds{1} \|} \leq 1
\end{equation}
where we used Theorem \ref{thm:numericalRoottest} to obtain the above inequality.  Thus $\alpha \leq \frac{\| 1 \|}{b}$ for arbitrary $b \in (0,R)$ and it follows that $\alpha \leq \frac{\| 1 \|}{R}$. $\Box$




\begin{exa} \label{exa:bandconverge}
Consider $f(z) = \sum_{n=0}^{ \infty} (1+j)z^n$ over the hyperbolic numbers $\Hcal = \RN \oplus j \RN$ where $j^2=1$. Observe $c_n= 1+j$ for all $n$ hence $\| c_n  \| = \sqrt{2}$ and the root test applies:
 \begin{equation}
 \alpha = \limsup_{n\to\infty}\sqrt[n]{\|c_n\|} = \lim_{ n \to \infty} 2^{1/2n} = 1.
 \end{equation}
Therefore, the Root Test provides for absolute convergence of this series where $\| z \| < 1/\Abound \alpha = 1 /\sqrt{2}$ as $\| zw \| \leq \sqrt{2} \| z \| \, \|w \|$ for hyperbolic numbers. Yet, consider $z = a(1-j)$ where $a \in \RN$,
 \begin{equation}
 f( a(1-j) ) = (1+j)(1)+ a(1+j)(1-j) + a^2(1+j)(1-j)^2+ \cdots = 1+j 
 \end{equation}
as $(1+j)(1-j)=0$. Therefore, $f(z)$ converges along the line $y=-x$ in the hyperbolic numbers where $z= x+jy$. Moreover, we can show this result extends to a whole band of nearby hyperbolic numbers. Consider points near $y=-x$ of the form $w=a+\varepsilon+j(\varepsilon-a)$. Such a point is distance $\sqrt{2} |\varepsilon|$ from $y=-x$. It is helpful to note the following identity:
 \begin{equation}
  (1+j)(x+yj) = (1+j)(x+y) \ \ \& \ \ (1+j)(x+yj)^n = (1+j)(x+y)^n 
 \end{equation}
for all $n \in \NN$. Thus, for $\varepsilon, a \in \RN$,
 \begin{equation}
 \|(1+j)w^n\| = \|(1+j)(a+\varepsilon + \varepsilon-a)^n\| = \sqrt{2}\, |2\varepsilon|^n.
\end{equation}
Using the root test for the numerical series $\sum (1+j)w^n$, we obtain
 \begin{equation}
 \alpha = \limsup_{n\to\infty}\sqrt[n]{\|c_n\|} = \lim_{ n \to \infty} |2\varepsilon| 2^{1/2n} = 2|\varepsilon|.
\end{equation}
This series converges for $|\varepsilon|<\frac12$ and thus $f(z$) converges for all $z$ at a distance of $\sqrt{2}|\varepsilon| < \frac1{\sqrt{2}}$ from the line $y=-x$.
Observe, this domain of convergence includes the disk $\|z\|<\frac1{\sqrt{2}}$ which we obtained via the root test, and, an infinite band of width $\sqrt{2}$ about the $y=-x$ line. 
\end{exa}

\noindent
Note, $\Abound = 1$ for $\Acal = \RN$ with $|x| = \sqrt{x^2}$ or $\Acal = \CN$ with $\| z \| = \sqrt{z \bar{z}}$ hence the result below is a natural extension of the usual geometric series from real or complex analysis; the radius of the domain for a geometric series in $\Acal$-calculus depends inversely on $\Abound$. 

\begin{thm} \label{thm:geometric}
For any commutative\footnote{this theorem is likely true in the noncommutative context as well} algebra $\Acal$, the geometric series $\sum z^n$ converges to $(1 - z)^{-1}$ for all $z\in\Acal$ such that $1-z\in\Acal^\times$ and $\|z\|< \frac1\Abound$.
\end{thm}

\noindent
\textbf{Proof:} The fact that the geometric series converges follows from the root test.  To show that it converges to $(1-z)^{-1}$, let $S_n = \sum_{k=0}^n z^k$ and $z$ meeting the conditions in the theorem.  Observe, 
\begin{align}
S_n(1-z) &= S_n-zS_n \\ \notag
&=1+z+z^2+...+z^n-z-z^2-...-z^n-z^{n+1}\\ \notag
&=1-z^{n+1}
\end{align}
and hence $S_n = (1-z^{n+1})(1 - z)^{-1}$ for all $n$ for which $(1-z)^{-1} \in \Acal^{\times}$. Thus,
\begin{align} 
\|S_n-(1 - z)^{-1}\| &= \|(1-z^{n+1})(1 - z)^{-1} - (1 - z)^{-1}\|\\ \notag
&=\|z^{n+1}(1 - z)^{-1}\|\\ \notag
&\leq \Abound\|z^{n+1}\|\|(1 - z)^{-1}\|\\ \notag
&\leq \Abound^{n+1}\|z\|^{n+1}\|(1 - z)^{-1}\| \ \ \ \text{(applied Proposition \ref{prop:ineqnpower})}
\end{align}
which can be made sufficiently small for $\|z\|<\frac1\Abound$.  $\square$ \\

\noindent
However, the domain defined in the theorem above imay fail to be maximal. Consider:

\begin{exa} \label{exa:directprodgeom}
Consider $\Acal  = \RN \times \RN$ where $(a,b) \star (x,y) = (ax,by)$ hence $(1,0) \star (0,1) = (0,0)$ are zero-divisors. Observe $\Acal^{ \times} = \{ (a,b) \ | \ ab \neq 0 \}$. Consider $z = (x,y)$ and 
\begin{equation} \sum_{n=0}^{ \infty} z^n = \sum_{n=0}^{ \infty}(x^n,y^n) = \left( \sum_{n=0}^{ \infty}x^n,\sum_{n=0}^{ \infty}y^n\right).
\end{equation}
The component series are geometric series with radii $x$ and $y$ respective. The series $\sum_{n=0}^{ \infty} z^n$ converges when both its component series converge. In particular, we need $|x|<1$ and $|y|<1$. We find the geometric series in $\RN \times \RN$ converges on the square $(-1,1)^2$. If we set $\| (x,y) \| = \sqrt{x^2+y^2}$ then it can be shown that $\| z \star w \| \leq \| z \| \|w \|$ hence $\Abound = 1$ for $\Acal = \RN \times \RN$ thus Theorem \ref{thm:geometric} only provides convergence on the disk $\| z \| <1$. 
\end{exa}

\noindent
In what follows we use the result of Example \ref{exa:bandconverge} to derive a result which is directly related to Example \ref{exa:directprodgeom} via the isomorphism of Proposition \ref{prop:isomorphismhyperbolic}.

\begin{exa}
In Example \ref{exa:bandconverge} we saw $f(z) =\sum_{n=0}^{\infty}(1+j)z^n$ converged on the infinite band 
\begin{equation}
 B_+ = \{ x+jy \ | \ -1-x \leq y \leq 1-x \}. 
\end{equation}
By entirely similar arguments, $g(z) = \sum_{n=0}^{\infty}(1-j)z^n$ will converge on
\begin{equation}
 B_- = \{ x+jy \ | \ -1+x \leq y \leq 1+x \}. 
\end{equation}
Note, $h(z) = f(z)+g(z)$ is defined for each $z \in B_+ \cap B_-$. In particular,
\begin{equation}
 h(z) = 2 \sum_{n=0}^{ \infty} z^n 
 \end{equation}
hence we find the geometric series in the hyperbolic numbers converges on a diamond with vertices $\pm 1, \pm j$. Notice, Theorem \ref{thm:geometric} merely provides convergence on the inscribed disk; $\| z \| < \frac{1}{\Abound} = \frac{1}{\sqrt{2}}$. Furthermore, the diamond with vertices $\pm 1$ and $\pm j$ is the image of the square $[-1,1]^2$ under the linear isomorphism $\Psi$ of Proposition \ref{prop:isomorphismhyperbolic}. 
\end{exa}

\noindent
The interplay between $\RN \times \RN$ with the direct product and the hyperbolic numbers $\RN \oplus j \RN$ is illustrative of an important calculational technique. It is often wise to exchange a problem in analysis in one algebra for a more lucid problem in an isomorphic algebra. \\

\noindent
Up to this point many of our theorems are likely true for noncommutative algebras. However, to discuss fractions in noncommutative algebras we would need to consider left and right divisors. We leave the noncommuting case to a future work.

\begin{thm} (\textbf{Ratio Test for Series with Unit-Coefficients }\label{thm:algebraratiotest}
Suppose $\Acal$ is an algebra with $ \| x \star y \| \leq \Abound \|x \| \|y \|$ for all $x,y \in \Acal$ and power series $\sum c_ n \star (z-z_o)^n$ where $c_n \in \Acal^\times$ for all $n$. Let
$$\alpha = \limsup_{n\to\infty}\left\|\frac{c_{n+1}}{c_n} \right\| \qquad \& \qquad R=\frac{1}{m_\Acal^2 \alpha}. $$
Then $\sum c_n \star (z-z_o)^n$ is absolutely convergent for $z-z_o \in \Acal^\times$ with $\|z-z_o\|<R$. Moreover, if $\alpha=0$ then $\sum c_n \star (z-z_o)^n$ converges absolutely on $\Acal$.

\end{thm}

\noindent
{\bf Proof:} Suppose $z_o=0$ and set $a_n=c_n \star z^n$ where $c_n,z \in \Acalx$ and note $a_n =c_n \star z^n \in \Acalx$ hence we are free to apply the ratio test for numerical series in $\Acal$: 
\begin{align}
\limsup_{n\to\infty}\frac{\|a_{n+1}\|}{\|a_n\|} &= 
\limsup_{n\to\infty}\frac{\|c_{n+1} \star z^{n+1}\|}{\|c_n \star z^n\|}  \\ \notag
&\leq \limsup_{n\to\infty}\Abound \left\|\frac{c_{n+1} \star z^{n+1}}{c_n \star z^n} \right\| \ \ \ \text{(by Corollary \ref{thm:quotientinequality})} \\ \notag
&=  \Abound\limsup_{n\to\infty} \left\|\frac{c_{n+1}}{c_n}\star z \right\| \\ \notag
&= \Abound\limsup_{n\to\infty}\Abound \left\|\frac{c_{n+1}}{c_n}  \right\| \| z \| \\ \notag
&= \Abound^2\| z \| \limsup_{n\to\infty} \left\|\frac{c_{n+1}}{c_n}  \right\|  \\ \notag
&= \Abound^2  \| z \| \alpha.
\end{align}
Since $R =\frac{1}{m_\Acal^2 \alpha}$ the condition $\Abound^2  \| z \| \alpha<1$ is equivalent to $\|z\| < R$ and Theorem \ref{thm:Ratio test} allows us to conclude that the series converges  absolutely for all $z \in \mathcal{A}^\times$ such that $\|z\| < R$. Finally, apply Lemma \ref{lem:centersub} to extend the proof to $z_o \neq 0$. $\Box$

\begin{exa}
Consider $\Hcal = \RN \oplus j \RN$ where $j^2=1$. Form the hyperbolic power series 
\begin{equation}
f(z) = \sum_{n=0}^{ \infty} (1+j)n! z^n.
\end{equation}
Observe,
\begin{equation}
 f( c(1-j)) = \sum_{n=0}^{ \infty} (1+j)n! c^n(1-j)^n = 1+j 
 \end{equation}
as $(1-j)(1+j)=0$ implies the terms with $n \geq 1$ all vanish. Thus the power series converges on $S=\{ c(1+j) \ | \ c \in \RN \}$. However, if $z \notin S$ then $(1+j)n! z^n \in \Hcal^{ \times}$ and
 it can be shown $\lim_{n \rightarrow \infty} (1+j)n! z^n \neq 0$ thus $f(z)$ diverges outside $S$. 
\end{exa}

\noindent
The Example above generalizes to other algebras. We can construct power series which converge on the zero-divisors or some subset of the zero-divisors and yet diverge everywhere else. Zero-divisors are simply beyond the scope of the ratio test for general power series in $\Acal$. Furthermore, Example \ref{exa:bandconverge} illustrates that the domain of convergence is not governed by root test alone. Interesting things can happen along zero divisors in the algebra. For example, we suspect $f(\zeta) = \sum_n (1+j+j^2)\zeta^n$ where $\zeta = x+yj+zj^2$ and $j^3=1$ converges on the infinite slab of thickness $2$ centered about the plane $x+y+z=0$. \\

\noindent
In contrast to Theorem \ref{thm:algebraratiotest}, the Theorem below can give us information about the convergence of the series at zero-divisors of the algebra.  

\begin{thm}(\textbf{Ratio Test for Series with Real Coefficients })\label{thm:algebraratiotestII}
Let $\Acal$ be an algebra with $ \| x \star y \| \leq \Abound \|x \| \|y \|$ for all $x,y \in \Acal$. power series $\sum c_n \star (z-z_o)^n$ where $0 \neq c_k \in \mathbb{R}$ for all $k \in \mathbb{N}$, put
$$\alpha = \limsup_{n\to\infty}\frac{|c_{n+1}|}{|c_n|},\ \ \ R=\frac{1}{m_\Acal\alpha}. $$
Then $\sum c_n \star (z-z_o)^n$ converges absolutely for $\|z-z_o\|<R$. Moreover, if $\alpha=0$ then $\sum c_n \star (z-z_o)^n$ converges absolutely on $\Acal$.
\end{thm}

\noindent
{\bf Proof:} Set $z_o=0$ to begin. If $c_n \in \RN$ then $ \| c_n \star z^n \| = |c_n| \| z^n \|$ . Put $a_n=c_n \star z^n$, and work towards applying the ratio test (Theorem \ref{thm:Ratio test}):
\begin{align}
\limsup_{n\to\infty}\frac{\|a_{n+1}\|}{\|a_n\|} &= \limsup_{n\to\infty}\frac{\|c_{n+1} z^{n+1}\|}{\|c_n z^n\|} \\ \notag
&= \limsup_{n\to\infty}\frac{|c_{n+1}|}{|c_n|} \frac{\| z^{n+1} \|}{\| z^n \|} \\ \notag
&\leq  \limsup_{n\to\infty}\frac{|c_{n+1}|}{|c_n|} \frac{ \Abound \|z\|\| z^{n} \|}{\| z^n \|} \\ \notag
&\leq  \Abound \|z\| \limsup_{n\to\infty}\frac{|c_{n+1}|}{|c_n|} \\ \notag
&= \Abound \|z\| \alpha. 
\end{align}
Since $R=\frac{1}{m_\Acal\alpha}$ we find $\Abound \|z\| \alpha < 1$ provides $ \| z \| < R$. Thus, applying Theorem \ref{thm:Ratio test}, the series converges for all $z \in \mathcal{A}^\times$ such that $\|z\| < R$. To conclude we apply Lemma \ref{lem:centersub} to extend the proof to $z_o \neq 0$. $\Box$ \\

\noindent
Next we compare and constrast the convergence indicated by Theorem \ref{thm:roottestpowerseries}, \ref{thm:algebraratiotest} and \ref{thm:algebraratiotestII}.

\begin{exa}
Consider $\sum_{n=1}^{\infty} \frac{3^n}{n} z^n$ for $z \in \Hcal = \RN \oplus j \RN$ with $j^2=1$. Notice, $c_n = \frac{3^n}{n}$ are real and recall $\Abound = \sqrt{2}$ for the hyperbolic numbers we consider here. Use $R_x$ to denote the radius of convergence suggested by Theorem $x$. We find:
\begin{equation}
 \left \| \frac{c_{n+1}}{c_n} \right \| = \frac{3n}{n+1} \rightarrow 3 \ \ \Rightarrow \ \ R_{\ref{thm:algebraratiotestII}} = \frac{1}{6} \ \ \& \ \  R_{\ref{thm:algebraratiotest}} = \frac{1}{3\sqrt{2}} 
\end{equation}
and,
\begin{equation}
\sqrt[n]{\| c_n \|} = \sqrt[n]{\frac{3^n}{n}} = \frac{3}{n^{1/n}} \rightarrow 3
\ \ \Rightarrow \ \ R_{\ref{thm:roottestpowerseries}} = \frac{1}{3\sqrt{2}}.
\end{equation}
Naturally, the root test Theorem \ref{thm:roottestpowerseries} and the real-coefficient ratio test Theorem \ref{thm:algebraratiotestII} both provide convergence of the series for $ \| z \| < \frac{1}{3 \sqrt{2}}$. However, Theorem \ref{thm:algebraratiotest} only indicates convergence for $ \| z \| < \frac{1}{6}$. There is nothing illogical about this as Theorem \ref{thm:algebraratiotest} is silent concerning the divergence of the series.
\end{exa}

\noindent
Next we study sequences and series of functions. In particular, we apply these results to power functions to gain insight into the theory of power series in $\Acal$. 

\begin{de}
Suppose $\{f_n\}$ is a sequence of functions $f_n: E \rightarrow \Acal$ where $E \subseteq \Acal$, and suppose the sequence of numbers $\{f_n(z)\}$ converges for every $z \in E$.  Then  
$f(z)=\lim_{n\to\infty}f_n(z) $
defines $f: E \rightarrow \Acal$. We say that $\{f_n\}$ converges to $f$ {\bf pointwise}.
\end{de}

\noindent
Pointwise convergence does not always guarantee properties of the sequence transfer to the limit. For example, it is possible to have  the limit function of a sequence of continuous functions which is discontinuous. To remedy this shortcoming of pointwise convergence we the stronger criteria of uniform convergence:

\begin{de}
We say that a sequence of functions $\{f_n\}$ on $E \subseteq \Acal$ {\bf converges uniformly to $f: E \rightarrow \Acal$} if for every $\varepsilon>0$ there is an integer $M$ such that $n\geq M$ implies
$\|f_n(z)-f(z)\| < \varepsilon$ for all $z\in E$. 
\end{de}

\noindent
Similar terminology is given for series of functions on $\Acal$.
We say that the {\bf series $\sum f_n(z)$  converges uniformly} on $E$ if the sequence $\{s_n\}$ of partial sums defined by
\begin{equation}
\sum_{i=1}^{n}f_i(z)=s_n(z)
\end{equation}
converges uniformly on E. There is also a Cauchy criterion for uniform convergence:
\begin{thm} \label{thm:cauchycriterionfnct}
A sequence of functions $\{f_n\}$ defined on $E \subseteq \Acal$ converges uniformly on $E$ if and only if for every $\varepsilon>0$ there exists $M \in \NN$ such that $m,n \geq M$ and $z\in E$ implies
$$\|f_n(z)-f_m(z)\| < \varepsilon.$$
\end{thm}

\noindent
{\bf Proof:} If $f_n \rightarrow f$ uniformly on $E$ then there exists an integer $M$ such that $n\geq M$ and  $z\in E$ imply
\begin{equation}
\|f_n(z)-f(z)\| < \frac{\varepsilon}{2}.
\end{equation}
Suppose $n, m \geq M$ and $z\in E$ and consider that
\begin{equation}
\|f_n(z)-f_m(z)\|\leq\|f_n(z)-f(z)\|+\|f_m(z)-f(z)\| < \varepsilon/2+\varepsilon/2 = \varepsilon. \end{equation}
Conversely, suppose the Cauchy condition holds.  As $\Acal$ is complete, the sequence $\{f_n(z)\}$ converges for every $z$, to a limit we may call $f(z)$.  Thus $f_n \rightarrow f$ on $E$.  It remains to show this convergence is uniform. Let $\varepsilon >0$ be given, and choose $M \in \NN$ such that $\|f_n(z)-f_m(z)\|< \varepsilon$.  Fix $n$ and let $m \rightarrow \infty$.  Since $f_m(z)\rightarrow f(z)$ as $m\rightarrow\infty$, this gives
\begin{equation}
 \|f_n(z)-f(z) \|  < \varepsilon
 \end{equation}
for every $n\geq M$ and  $z\in E$ hence $f_n \rightarrow f$ uniformly on $E$. $\square$\\ 

\noindent
Weierstrauss' taught us that uniform convergence of series of functions on $E \subseteq \CN$ can be derived from the existence of a {\it majorizing series}. In particular, if a convergent numerical series bounds the values of the function on $E$ then the series of functions converges uniformly on $E$. This is often known as Weierstrauss $M$-test. 

\begin{thm} \label{thm:weierstraussM}
( {\bf Weierstrauss $M$-Test for $\Acal$} ) Suppose $\{f_n(z)\}$ is a sequence of functions defined on $E$. If $\sum M_n$ is a convergent series in $\RN$ and $ \|f_n(z) \|\leq M_n$ for all $z\in E$ and $n\in\NN$ then $\sum f_n$ converges uniformly on $E$.
\end{thm}

\noindent
{\bf Proof:} Assume  $ \|f_k(z) \|\leq M_k$ for each $k \in \NN$ and note for $m \geq n$: 
\begin{equation}
  \left\|\sum_{i=1}^{m}f_i(z) - \sum_{i=1}^{n}f_i(z) \right\| = \left\|\sum_{i=n}^{m}f_i(z) \right\| \leq \sum_{i=n}^{m} \left\|f_i(z) \right\| \leq \sum_{i=n}^{m} M_n
  \end{equation}
for each $z \in E$. Furthermore, by convergence of $\sum M_n$, for each $\varepsilon>0$ we may select $M \in \NN$ for which $m,n \geq M$ imply $\sum_{i=n}^{m} M_n < \varepsilon$. Consequently, the conditions of Theorem \ref{thm:cauchycriterionfnct} are met and we conclude $\sum f_n$ converges uniformly on $E$. $\square$

\begin{thm} \label{thm:pwrseriesnormalconvergence}
If $\sum c_n \star (z-z_o)^n$ is a power series on $\Acal$ and
$$\alpha = \limsup_{n\to\infty}\sqrt[n]{\|c_n\|} \qquad \& \qquad  \ R=\frac{1}{ \Abound \alpha }. $$ 
then $\sum_n c_n \star (z-z_o)^n $ is uniformly absolutely convergent for $\|z-z_o\| \leq  R- \varepsilon$ for any $\varepsilon \in (0, R)$.  If R is infinite, then the series is uniformly absolutely convergent for $\|z\|\leq  L$ for any $L>0$.
\end{thm}

\noindent
{\bf Proof:} absolute convergence is given by Theorem \ref{thm:roottestpowerseries}. Let $\varepsilon \in (0, R)$ and choose $z_1$ such that $z_1-z_o = \Abound(R- \varepsilon) > 0$. Note $z_1-z_o$ is by construction real and $z_1-z_o < \Abound R = 1/\alpha$ hence: 
\begin{equation}
 \limsup_{n\to\infty}\sqrt[n]{\|c_n \star (z_1-z_o)^n\|} = (z_1-z_o)\limsup_{n\to\infty}\sqrt[n]{\|c_n\|} < \alpha /\alpha  = 1.
\end{equation}
Therefore, $\sum_n \| c_n \star (z_1-z_o)^n \|$ converges by Theorem \ref{thm:numericalRoottest}. For, $z$ with $\| z-z_o \| \leq R- \varepsilon$,
 \begin{align}  
 \| c_n \star (z-z_o)^n \| &\leq \Abound^n \|c_n\| \|z-z_o\|^n \\ \notag
 &\leq \|c_n \| \Abound^n(R- \varepsilon)^n \\ \notag
 &= \|c_n \| (z_1-z_o)^n \\ \notag
 &=\|c_n\star (z_1-z_o)^n\|
 \end{align}
 for each $n$ thus by the Weierstrauss $M$-Test (\ref{thm:weierstraussM}) the theorem follows $\Box$. \\
    
\noindent
Integration in $\Acal$ and uniform limits can be interchanged:

\begin{thm} \label{thm:uniformcontinuityintegral} 
Suppose $C$ be a piecewise smooth curve of length $L< \infty$ in $\Acal$ and suppose $U$ is an open set containing $C$. If $\{ f_j \}$ is a sequence of continuous $\Acal$-valued functions on $U$, and if $\{ f_j \}$ converges uniformly to $f$ on $U$ then $\int_{C} f_j(z) \star dz$ converges to $\int_{C} f(z) \star dz$.
\end{thm}

\color{black}
\noindent
{\bf Proof:} suppose $\epsilon >0$. Note, uniform convergence of $\{ f_j \}$ to $f$ implies there exists $N \in \NN$ for which $n \geq N$ implies $ ||f_n(z)-f(z)|| < \frac{\epsilon}{\Abound L} $ for all $z \in U$. Since $C \subset U$ we have the same estimate for points on $C$. Furthermore, by the treatment of integration in \cite{cookAcalculusI},
\begin{align}
\bigg{|}\bigg{|} \int_C f_n(z) \star dz - \int_C f(z) \star dz\bigg{|}\bigg{|} &= \bigg{|}\bigg{|} \int_C \left( f_n(z)-f(z) \right) \star dz\bigg{|}\bigg{|} \\ \notag
&\leq \Abound\frac{\epsilon}{\Abound L}L = \epsilon. \ \ \Box 
\end{align}

\noindent
Observe: if $f_n \rightarrow f$ uniformly near $C$ then $\lim_{n \rightarrow \infty} \int_C f_n \star dz = \int_C \left( \lim_{n \rightarrow \infty} f_n \right) \star dz$. \\

\noindent
In complex analysis we learn one consquence of Cauchy's Integral Formula is that a uniformly convergent sequence of complex-differentiable functions has a limit function which is likewise complex-differentiable. However, in the absense of Cauchy's Integral Formula, no such luxury is available in the study of differentiability of the limit function. We face the usual difficulty of real analysis which is nicely addressed by Dieudonn\a' e in result 8.6.3 of \cite{Dmaster}. 
We show Dieudonn\a' e's result extends naturally to $\Acal$-calculus: we provide sufficient conditions for a sequence of $\Acal$-differentiable functions to have an $\Acal$-differentiable limit function:

\begin{thm} \label{thm:difflimitfunseries}
Let $U$ be an open connected subset of $\Acal$, $f_n: U \rightarrow \Acal$ an $\Acal$-differentiable mapping of $U$ for each $n \in \NN$. Suppose that: 
\begin{enumerate}[{\bf (i.)}]
\item there exists one point $z_0 \in U$ such that the sequence $\{ f_n(z_0) \}$ converges in $\Acal$,
\item for every point $a \in U,$ there is a ball $B(a)$ of center $a$ contained in $U$ and such that in $B(a)$ the sequence $\{ f'_n \}$ converges uniformly. 
\end{enumerate}
Then for each $a \in U,$ the sequence $\{f_n \}$ converges uniformly in $B(a)$; moreover, if, for each $z \in U,$ $f(z)=\lim_{n \rightarrow \infty} f_n(z)$ and $g(z) = \lim_{n \rightarrow \infty} f'_n(z),$ then $g(z)=f'(z),$ for each $z \in U.$ 
\end{thm}

\noindent
To be clear, when we write $f'(z)$ this indicates the $\Acal$-derivative of $f$. Hence, in part, the Theorem asserts $f_n \rightarrow f$ where $f$ is $\Acal$ differentiable. Furthermore, for each $z \in U$:
$$ \frac{d}{dz} \left(\lim_{n \rightarrow \infty}f_n(z) \right) =  \lim_{n \rightarrow \infty} \left( \frac{df_n}{dz}(z) \right). $$

\noindent
{\bf Proof:} suppose $f_n: U \rightarrow \Acal$ is an $\Acal$-differentiable mapping on an open connected $U \subseteq \Acal$. In addition, suppose conditions (i.) and (ii.)  hold. Notice that $\Acal$-differentiable implies Frechet differentiable. Hence, by result 8.6.3 in  \cite{Dmaster} we find uniform convergence of $\{ f_n \}$ as described in the Theorem. Dieudonn\a' e uses the  notation $f'(x)$ for the Frechet derivative of $f$ at $x$. We change notation and write $Df(x)$ for the Frechet derivative of $f$ at $x$. Hence, by $(8.6.3)$ in \cite{Dmaster}, if for each $z \in U,$ $f(z)=\lim_{n \rightarrow \infty} f_n(z)$ and $g(z) = \lim_{n \rightarrow \infty} Df_n(z),$ then $g(z)=Df(z),$ for each $z \in U$. Let $\{ v_1, v_2, \dots , v_N \}$ serve as a basis for $\Acal$ with coordinates $x_1, x_2, \dots , x_n$. By the definition of partial derivative, for each $i=1,2, \dots , N$,
\begin{equation}
Df(z)(v_i) = \frac{\partial f}{\partial x_i}(z) \qquad \& \qquad Df_n(z)(v_i) = \frac{\partial f_n}{\partial x_i}(z) 
\end{equation}
thus \cite{Dmaster} provides the existence of the Frechet derivative as well as the following identity for the partial derivatives:
\begin{equation} \label{eqn:exchangelimits}
\frac{\partial }{\partial x_i} \left(  \lim_{n \rightarrow \infty} f_n(z)  \right) = \lim_{n \rightarrow \infty}  \left( \frac{\partial f_n}{\partial x_i}(z) \right).
\end{equation}
It remains to show $f = \lim_{n \rightarrow \infty} f_n$ is $\Acal$-differentiable on $U$.  From (ii.) we know $f_n$ is $\Acal$-differentiable hence $f_n$ satisfy the symmetric CR-equations\footnote{if $v_1 = \mathds{1}$ then symmetric CR-equations reduce to the usual CR-equations $\frac{\partial f}{\partial x_j} = \frac{\partial f}{\partial x_1} \star v_j$. } 
\begin{equation}
 \frac{\partial f_n}{\partial x_i} \star v_j = \frac{\partial f_n}{\partial x_j} \star v_i.
 \end{equation}
Hence, using the symmetric $\Acal$-CR equations and Equation \ref{eqn:exchangelimits} we derive:
\begin{align}
 \frac{\partial f}{\partial x_i} \star v_j &= \frac{\partial}{\partial x_i}\left[\lim_{n \rightarrow \infty} f_n \right] \star v_j  \\ \notag
 &= \lim_{n \rightarrow \infty} \left[ \frac{\partial f_n}{\partial x_i} \right]  \star v_j \\ \notag
 &= \lim_{n \rightarrow \infty} \left[ \frac{\partial f_n}{\partial x_i}   \star v_j \right] \\ \notag
 &=  \lim_{n \rightarrow \infty} \left[ \frac{\partial f_n}{\partial x_j}   \star v_i \right].
 \end{align}
Consequently, $\frac{\partial f}{\partial x_i} \star v_j = \frac{\partial f}{\partial x_j} \star v_i$ and thus $f$ is $\Acal$-differentiable with $g(z) = f'(z)$. $\Box$ \\

\noindent
Theorem \ref{thm:difflimitfunseries} allows us to establish the $\Acal$-differentiability of power series in $\Acal$. In particular, if an $\Acal$-series converges on an open ball about its center then we find the derivative of the series exists and can be obtained by term-wise differentiation.

\begin{coro} \label{coro:derivativepowseries}
If $\sum c_n \star (z-z_o)^n$ is a power series on $\Acal$ which converges for $\| z - z_o \| < R$. Then $ \frac{d}{dz} \sum c_n \star (z-z_o)^n = \sum nc_n \star (z-z_o)^{n-1}$ for each $z \in \Acal$ with $\| z-z_o \| < \frac{1}{\Abound \alpha}$ where $\alpha = \limsup_{n\to\infty}\sqrt[n]{\|c_n\|}$. 
\end{coro}

\noindent
{\bf Proof:} We begin with $z_o =0$. Assume the series $\sum c_n \star z^n$ converges for $\| z \| < R$.  Let $U = \{ z \in \Acal \ | \ \| z \| <R \}$ and note $U$ is an open set. Define $f_n(z) = \sum_{k=0}^n c_k \star z^k$ for $n=0,1,\dots $ and $z \in U$. Observe $\frac{d f_n}{dz} = \sum_{k=1}^n kc_k \star z^{k-1}$ for $z  \in U$ and $0 \in U$ with $\{ f_n (0) \} = \{ c_0 \}$ is convergent. If $a \in U$ then note $B_{r}(a) = \{ z \ | \ \| z-a \| <r \} \subseteq U$ where $r = \text{min}\{ \| a \|/2, |R-\| a \||/2 \}$. We need to show $\{ f_n' \}$ converges uniformly on $B_r(a)$. Let $\alpha =\limsup_{n\to\infty}\sqrt[n]{\|c_n\|}$ then Theorem \ref{thm:norootcriterianeeded} provides $\alpha \leq \frac{\mathds{1}}{R}$ which shows $\alpha$ is finite. Since $\sqrt[n]{n}\rightarrow 1\ as\ n\to\infty$, we have 
\begin{equation}
\limsup_{n\to\infty}\sqrt[n]{n\|c_n\|}=\limsup_{n\to\infty}\sqrt[n]{\|c_n\|} = \alpha < \infty
\end{equation}
hence Theorem \ref{thm:roottestpowerseries} provides that $\sum nc_n \star z^{n-1}$ converges for $\| z \| <  \frac{1}{\Abound \alpha}$. Therefore, by Theorem \ref{thm:pwrseriesnormalconvergence}, we find $\{ f_n' \}$ is uniformly convergent for $\| z \| < R- \varepsilon$ for any $\varepsilon \in (0, R)$. Thus $\{ f_n' \}$ converges uniformly on $B_r(a)$ as we are free to adjust $\varepsilon$ such that $B_r(a) \subset B_{R-\epsilon}(0)$. In summary, we have satisfied conditions (i.) and (ii.) of Theorem \ref{thm:difflimitfunseries} and the Corollary follows from the identifications $f(z) = \sum c_n \star z^n$ and $g(z) = \sum nc_n \star z^{n-1}$ for which $f'(z)=g(z)$ on $U$. Finally, we apply Lemma \ref{lem:centersub} to shift the $z_o=0$ result to $z_o \neq 0$. $\Box$ \\

\noindent
A similar result is available for higher derivatives:

\begin{coro} \label{coro:coeffbyderivatives}
If $f(z) =\sum c_n \star (z-z_o)^n$ is a power series on $\Acal$ which converges for $\| z-z_o \| < R$ then $f$ has derivatives of all orders in $\|z-z_o\|<\frac{1}{\Abound \alpha}$ where $\alpha = \limsup_{n\to\infty}\sqrt[n]{\|c_n\|}$. Moreover, the higher-derivative functions are obtained by term-wise differentiation:
$$f^{(k)}(z)=\sum_{n=k}^\infty n(n-1)\cdot \cdot \cdot (n-k+1)c_n \star (z-z_o)^{n-k}$$
and $f^{(k)}(z_o)=k!c_k$ for $k=0, 1, 2, \dots$.
\end{coro}

\noindent
{\bf Proof:} observe $
\lim_{n \rightarrow \infty}\sqrt[n]{n(n-1)\cdot \cdot \cdot (n-k+1)} =1$ for any $k \in \NN$. Thus the argument given for Corollary \ref{coro:derivativepowseries} naturally extends to the $k$-th derivative. $\square$ \\

\noindent
We close this section with a discussion of entire functions on an $\Acal$.

\begin{de}
A function $f:\Acal\to\Acal$ is called entire if it can be written as a power series $\sum a_n \star z^n$ which converges on all of $\Acal$.
\end{de}

\begin{thm}\label{thm:entire}
If $f(x) = \sum a_n x^n$ is an entire function on the reals, then there exists a unique entire extension to $\Acal$. The extension has the form $\tilde{f}(z) = \sum a_n z^n$ for each $z \in \Acal$.
\end{thm}

\noindent
\textbf{Proof:} Let $f$ be an entire function on the reals.  Thus, its radius of convergence is infinite, and by the real root test we have:
\begin{equation}
\limsup_{n\to\infty} \sqrt[n]{|a_n|} = 0
\end{equation}
Let $\tilde{f}(z) = \sum a_n z^n$ for each $z \in \Acal$.  Since the coefficients of the extended function $\tilde{f}(z)$ are the same as those of $f(x)$ we find $\tilde{f}(z)$ is entire by Theorem \ref{thm:roottestpowerseries}. If $g(z) = \sum b_n \star z^n$ is entire function on $\Acal$ for which $g|_{\RN} = f= \tilde{f}|_{\RN}$ then $g^{(n)}(0) = \tilde{f}^{(n)}(0)$ for $n=0,1,2,\dots $. Hence, $ b_n = a_n$ for $n=0,1,2,\dots $ by Corollary \ref{coro:coeffbyderivatives}. Thus the extension $\tilde{f}(z) = \sum a_nz^n$ is the unique extension to $\Acal$. $\square$ \\



\begin{thm} \label{thm:entirefnctuniformabs}
If f is an entire function on $\Acal$, then we have
$$\limsup_{n\to\infty}\sqrt[n]{\|c_n\|}=0$$
and thus f is uniformly absolutely convergent for $\|z\|<L$ for all $L>0$.
\end{thm}

\noindent
\textbf{Proof:} Since $f$ is entire, it follows that $\sum c_n z^n$ converges on $\Acal$ and hence, by Theorem \ref{thm:term test}, we have 
\begin{equation}
\lim_{n\to\infty} c_n z^n = 0 \Rightarrow \lim_{n\to\infty} \|c_n z^n\| = 0 \qquad \forall z\in\Acal.
\end{equation}
Given $\varepsilon>0$, there exist $N\in\NN$ such that for each $n\geq N$,
\begin{equation}
\left\|c_n\left(\frac1\varepsilon\right)^n \right\|< 1 \ \ 
\Rightarrow \ \ \|c_n\|<\varepsilon^n \ \ \Rightarrow \ \  \sqrt[n]{\|c_n\|}<\varepsilon. 
\end{equation}
Thus, $\limsup_{n\to\infty}\sqrt[n]{\|c_n\|}=0$ and by Theorem \ref{thm:pwrseriesnormalconvergence} we reach the desired result.$\ \square$

\begin{coro}
The set of entire functions on $\Acal$ is an algebra, and the product of two entire functions $\sum a_n \star z^n$ and $\sum b_n \star z^n$ is equal to $\sum c_n \star z^n$ where $c_n=\sum_{k=0}^{n}a_k \star b_{n-k}$.
\end{coro}

\noindent
\textbf{Proof:} It is clear that this set is a real vector space, so all we must show is that the product of two entire functions $\sum a_n \star z_n$ and $\sum b_n \star z_n$ is also entire.  By Theorem \ref{thm:pwrseriesnormalconvergence}, these functions are absolutely convergent on all of $\Acal$, and thus by Theorem \ref{thm:productofseries}, their product will converge on all of $\Acal$ and will be  
\begin{align}
\left(\sum a_n \star z^n\right) \star \left(\sum b_n  \star z^n\right)&= \sum_{n=0}^\infty\left(\sum_{k=0}^n a_k \star z^k\star b_{n-k} \star z^{n-k}\right)\\ \notag
&=\sum_{n=0}^\infty\left(\sum_{k=0}^n a_k \star b_{n-k} \star z^n\right)\\ \notag
&=\sum_{n=0}^\infty\left(\sum_{k=0}^n a_k \star b_{n-k}\right) \star z^n\\ \notag
&=\sum_{n=0}^\infty c_n\star z^n \  where\ c_n=\sum_{k=0}^{n}a_k  \star b_{n-k}.
\end{align}
Thus, we have the desired result.  $\square$

\section{Transcendental functions} \label{sec:transfunct}
In this section we let $\Acal$ denote a real associative finite dimensional commutative algebra $\Acal$. Theorem \ref{thm:entire} encourages us to take the power series formulation of elementary functions as fundamental. If we want to recover standard elementary functions by restriction to $\RN \subseteq \Acal$ then our definitions must be given. We provide concrete definitions for the exponential, sine, cosine, hyperbolic sine and hyperbolic cosine over any $\Acal$. We also provide proofs for identities of these functions which equally well apply to a myriad of choices for $\Acal$. 

\subsection{Exponential}
\begin{de}
For each z $\in\Acal$, we define
$$exp(z)=\sum_{n=0}^\infty \frac{z^n}{n!}.$$
\end{de}
\noindent

\begin{thm} \label{thm:propofexp}
Let $\Acal$ be a commutative, unital, associative algebra over $\RN$.   
\begin{quote}
\begin{enumerate}[{\bf (i.)}]
\item $exp(z) = \sum_{n=0}^\infty \frac{z^n}{n!}$ is entire on $\Acal$,
\item $exp(z)\star exp(w) = exp(z+w) $ for all $z,w \in \Acal$,
\item $exp(0)=1$ and $exp(-z) =exp(z)^{-1}$ hence $exp(z) \in \Acalx$.
\end{enumerate}
\end{quote}
\end{thm}

\noindent
{\bf Proof:} {\bf (i.)} Identify the coefficients of the exponential are $c_n=\frac{1}{n!}$, which gives us
\begin{equation}
\limsup_{n\to\infty}\sqrt[n]{\|c_n\|}=\limsup_{n\to\infty}\frac1{\sqrt[n]{n!}}=0.
\end{equation}
Thus, by Theorem \ref{thm:pwrseriesnormalconvergence}, this series uniformly converges absolutely for all $z\in\Acal$. \\

\noindent
 {\bf (ii.)} By Theorem \ref{thm:productofseries}, we have for each $z, w \in\Acal$:
\begin{align}
exp(z)\star exp(w)&=\biggl{(}\sum_{n=0}^\infty \frac{z^n}{n!}\biggr{)}\star \biggl{(}\sum_{n=0}^\infty \frac{w^n}{n!}\biggr{)}\\ \notag
&=\sum_{n=0}^\infty \sum_{k=0}^n \frac{z^k}{k!}\star\frac{w^{n-k}}{(n-k)!}\\ \notag
&=\sum_{n=0}^\infty \sum_{k=0}^n \frac{n!}{k!(n-k)!}\frac{z^n\star w^{n-k}}{n!}\\ \notag
&=\sum_{n=0}^\infty \frac{1}{n!}\sum_{k=0}^n {n\choose k}z^n\star w^{n-k}\\ \notag
&=\sum_{n=0}^\infty\frac{(z+w)^n}{n!}\\ \notag
&=exp(z+w).
\end{align}
{\bf (iii.)} Clearly $exp(0)=1$ and since for each $z\in\Acal$ we have 
\begin{equation}
exp(z) \star exp(-z)=exp(z-z)=exp(0)=1, 
\end{equation}
thus $exp(z)\in\Acalx$ for each $z \in \Acal$. Moreover,
$exp(-z)=exp(z)^{-1}$. $\Box$ \\

\noindent
The following Theorem should not be surprising.

\begin{thm}
For each $z \in \Acal$, $\frac{d}{dz}exp(z) = exp(z)$.
\end{thm}

\noindent
{\bf Proof:} From Corollary \ref{coro:derivativepowseries}, we know that $\frac{d}{dz}exp(z)$ exists and
\begin{equation}
\frac{d}{dz}exp(z)=\sum_{n=1}^\infty \frac{nz^{n-1}}{n!}=\sum_{n=0}^\infty \frac{z^{n}}{n!}=exp(z)
\end{equation}
for all $z\in\Acal$. $\Box$ \\

\noindent
To appreciate the content of this section we expand the exponential into component functions of several interesting algebras.
\begin{exa}
Let $\Acal = \RN$ then $exp(x)=e^x$ is the usual exponential of real calculus. 
\end{exa}

\begin{exa}
Let $\Acal = \CN$ then $exp(x+iy)= exp(x)exp(iy) = e^x \cos y+ ie^x \sin y$. Setting $u+iv = exp(x+iy)$ we find $u= e^x \cos y$ and $v=e^x \sin y$. The component functions of the exponential are solutions to both the Cauchy Riemann equations $u_x=v_y$ and $v_x = -u_y$ and Laplace's Equation $\phi_{xx}+\phi_{yy}=0$. 
\end{exa}

\begin{exa}
Let $\Acal = \Hcal$ then $exp(x+jy)= exp(x)exp(jy) = e^x( \cosh y+ j \sinh y)$. Setting $u+jv = exp(x+jy)$ provides $u = e^x \cosh y$ and $v=e^x\sinh y$. These functions satisfy the hyperbolic Cauchy Riemann equations $u_x=v_y$ and $u_y = v_x$ and both are solutions to the wave equation $\phi_{xx}-\phi_{yy}=0$.
\end{exa}

\begin{exa}
Let $\Acal = \RN \oplus \eps \RN$ with $\eps^2=0$. Then $exp(x+\eps y)= exp(x)exp(\eps y)$. Calculate
\begin{equation}
 exp( \eps y) = 1+ \eps y+ \frac{1}{2}\eps^2 y^2+ \cdots =  1+ \eps y 
\end{equation}
hence $exp(x+\eps y) = e^x+ \eps ye^x $. The component functions of the exponential are $e^x$ and $ye^x$ for this nilpotent algebra. \end{exa}

\subsection{Hyperbolic sine and cosine}
As in the usual calculus, the hyperbolic sine and cosine appear as the odd and even pieces of the exponential function.

\begin{de} \label{def:hyperbolictrig}
For each $z \in \Acal$, we define
$$ cosh(z)=\sum_{n=0}^\infty \frac{z^{2n}}{(2n)!} \qquad \& \qquad sinh(z)=\sum_{n=0}^\infty \frac{z^{2n+1}}{(2n+1)!}. $$
\end{de}

\begin{thm} \label{thm:convhyperbolic}
The series defining $cosh(z)$ and $sinh(z)$ uniformly converge absolutely on $\Acal$. \end{thm}

\noindent
{\bf Proof:} For $cosh(z)$ we have coefficents $c_n = 1/(2n)!$ hence calculate
\begin{equation}
\limsup_{n\to\infty}\sqrt[n]{\|c_n\|}=\limsup_{n\to\infty}\frac1{\sqrt[n]{(2n)!}}=0
\end{equation}
and for $sinh(z)$ we have coefficients $c_n = 1/(2n+1)!$ hence calculate:
\begin{equation}
\limsup_{n\to\infty}\sqrt[n]{\|c_n\|}=\limsup_{n\to\infty}\frac1{\sqrt[n]{(2n+1)!}}=0
\end{equation}
Thus, by the Theorem \ref{thm:pwrseriesnormalconvergence} we find both series converge absolutely for each $ z \in \Acal$. $\Box$

\begin{thm} \label{thm:hyperbolicpropsI}
Hyperbolic sine and cosine have the following properties:
\begin{quote}
\begin{enumerate}[{\bf (i.)}]
\item $cosh(0)=1$ and $sinh(0)=0$,
\item $cosh(-z) = cosh(z)$ and $sinh(-z) = -sinh(z)$ for each $z \in \Acal$,
\item $\frac{d}{dz} cosh(z) = sinh(z)$ and $\frac{d}{dz} sinh(z) = cosh(z)$ for each $z \in \Acal$.
\end{enumerate}
\end{quote}
\end{thm}

\noindent
{\bf Proof:} Observe {\bf (i.) } and {\bf (ii.)} follow immediately from Definition \ref{def:hyperbolictrig}. Using Corollary \ref{coro:derivativepowseries} we derive {\bf (iii.)} as follows:
\begin{align}
\frac{d}{dz}cosh(z)&=\sum_{n=1}^\infty \frac{2nz^{2n-1}}{(2n)!}=\sum_{n=1}^\infty \frac{z^{2n-1}}{(2n-1)!}=\sum_{n=0}^\infty \frac{z^{2n+1}}{(2n+1)!}=sinh(z) \\ \notag
\frac{d}{dz}sinh(z)&=\sum_{n=0}^\infty \frac{(2n+1)z^{2n}}{(2n+1)!}=\sum_{n=0}^\infty \frac{z^{2n}}{(2n)!}=cosh(z). \ \ \Box 
\end{align}

\noindent
Sums and products of these series will also be entire, giving us the following:
\begin{thm} \label{thm:hyperbolicpythago}
$cosh^2(z)-sinh^2(z)=1$ for all $z\in\Acal$.
\end{thm}

\noindent
{\bf Proof:} Let $g(z)=cosh^2(z)-sinh^2(z)$ for all $z\in\Acal$.  By Theorems \ref{thm:hyperbolicpropsI}, \ref{thm:convhyperbolic} and the chain-rule, 
\begin{equation}
g'(z)=2cosh(z) \star sinh(z)-2sinh(z) \star  cosh(z)=0 
\end{equation}
for all $z \in\Acal$. Since $\Acal$ is a connected  it follows $g(z)$ is constant. Moreover, 
\begin{equation}
 g(0)=cosh^2(0)-sinh^2(0)=1
\end{equation}
hence the Theorem follows.  $\square$

\begin{thm} \label{thm:hyperbolicidentI}
For all $z\in\Acal$,
\begin{quote}
\begin{enumerate}[{\bf (i.)}]
\item  $e^z=cosh(z)+sinh(z)$,
\item  $cosh(z)=\frac{1}{2}(e^z+e^{-z})$,
\item  $sinh(z)=\frac{1}{2}(e^z-e^{-z})$.
\end{enumerate}
\end{quote}
\end{thm}

\noindent
{\bf Proof:} item (i.) is verified directly from the definitions: 
\begin{equation}
cosh(z)+sinh(z)=\sum_{n=0}^\infty \frac{z^{2n}}{(2n)!}+\sum_{n=0}^\infty \frac{z^{2n+1}}{(2n+1)!}=\sum_{n=0}^\infty \frac{z^{n}}{n!}=e^z.
\end{equation}
Hence (ii.) follows from (i.) by simple calculation:
\begin{align}
\frac{1}{2}(e^z+e^{-z}) &= \frac{1}{2}(cosh(z)+sinh(z)+cosh(-z)+sinh(-z))\\ \notag
&= \frac{1}{2}(cosh(z)+sinh(z)+cosh(z)-sinh(z)) \\ \notag
&= cosh(z).
\end{align}
Likewise (iii.) follows from (i.)
\begin{align}
\frac{1}{2}(e^z-e^{-z})&= \frac{1}{2}(cosh(z)+sinh(z)-cosh(-z)-sinh(-z))\\ \notag
&= \frac{1}{2}(cosh(z)+sinh(z)-cosh(z)+sinh(z)) \\ \notag
&= sinh(z).
\end{align}
Alternatively, we could have differentiated (ii.) to obtain (iii.). $\Box$

\begin{thm}
For all $z, w\in\Acal$
$${\bf (i.)} \ \ cosh(z+w)=cosh(z)\star cosh(w)+sinh(a)\star sinh(w),$$
$${\bf (ii.)} \ \ sinh(z+w)=sinh(z)\star cosh(w)+cosh(z)\star sinh(w).$$
\end{thm}

\noindent
{\bf Proof:} For all $z, w\in\Acal$ apply Theorems \ref{thm:hyperbolicidentI} and \ref{thm:propofexp} part (ii.)
\begin{align}
cosh(z+w)&=\frac{1}{2}(e^{z+w}+e^{-(z+w)})\\ \notag
&=\frac{1}{2}(e^z\star e^w+e^{-z}\star e^{-w})\\ \notag
&=\frac{1}{2}(e^z \star e^w+e^z \star e^{-w}-e^z \star e^{-w}+e^{-z}\star e^{-w})\\
\notag
&=\frac{1}{2}(e^z\star (e^w+e^{-w})-(e^z-e^{-z}) \star e^{-w})\\
\notag
&=\frac{1}{2}(2e^z \star cosh(w)-2sinh(z) \star e^{-w})\\
\notag
&= c(z) \star c(w)+s(z) \star c(w)-s(z) \star c(-w)-s(z) \star s(-w)\\
\notag
&=cosh(z)\star cosh(w)+sinh(z)\star sinh(w)
\end{align}
where we used the abbreviated notation $c(z) =cosh(z)$ and $s(z)=sinh(z)$ in the next to last line. Fix $w$ and differentiate with respect to $z$ to find
$$sinh(z+w)=sinh(z) \star cosh(w)+cosh(z) \star sinh(w). $$
Thus the adding angles formulas for hyperbolic functions exist for $\Acal$. $\Box$

\subsection{Sine and cosine}
\noindent
In certain contexts we could use the imaginary unit $i$ for which $i^2=-1$ to aid in the definition of sine and cosine. However, there are many algebras without such an imaginary unit hence we provide a treatment which only requires the theory of power series to establish the structure of trigonometric functions. Once more, we find identities for trigonometric function which transcend the choice of $\Acal$. Our arguments here are parallel those found in Section 68 of \cite{JandP}.

\begin{de} \label{def:sinecosine}
For each z $\in\Acal$, we define
$$cos(z)=\sum_{n=0}^\infty (-1)^n\frac{z^{2n}}{(2n)!}
\qquad \& \qquad sin(z)=\sum_{n=0}^\infty (-1)^n\frac{z^{2n+1}}{(2n+1)!}.$$
\end{de}

\begin{thm} \label{thm:convtrig}
The series defining $cos(z)$ and $sin(z)$ converge absolutely on $\Acal$. \end{thm}

\noindent
{\bf Proof:} For $cos(z)$ note $c_n = (-1)^n/(2n)!$ and for $sin(z)$ the coefficient $c_n = (-1)^n /(2n+1)!$ consequently the proof for Theorem \ref{thm:convhyperbolic} equally well applies here. $\Box$  

\begin{thm} \label{thm:trigpropsI}
Sine and cosine over $\Acal$ have the following properties:
\begin{quote}
\begin{enumerate}[{\bf (i.)}]
\item $cos(0)=1$ and $sin(0)=0$,
\item $cos(-z) = cos(z)$ and $sin(-z) = -sin(z)$ for each $z \in \Acal$,
\item $\frac{d}{dz} cos(z) = -sin(z)$ and $\frac{d}{dz} sin(z) = cos(z)$ for each $z \in \Acal$.
\end{enumerate}
\end{quote}
\end{thm}

\noindent
{\bf Proof:} follows from argument nearly identical to those given for Theorem \ref{thm:hyperbolicpropsI}. $\Box$ \\

\noindent
The following identity is well-known for $\RN$ or $\CN$, but it just as well applies to sine and cosine over any $\Acal$:

\begin{thm} \label{thm:pythagcosinesine}
$cos^2(z)+sin^2(z)=1$ for all $z\in\Acal$.
\end{thm}

\noindent
{\bf Proof:} Let $g(z)=cos^2(z)+sin^2(z)$ for all $z\in\Acal$.  By Theorems \ref{thm:convtrig} and \ref{thm:trigpropsI} and the chain-rule, we have for all $z \in\Acal$
\begin{equation}
g'(z)=-2cos(z)\star sin(z)+2sin(z)\star cos(z)=0. 
\end{equation}
Thus $g(z)$ is constant on all of $\Acal$.  Since $g(0)=cos^2(0)+sin^2(0)=1$, the result follows. $\square$

\begin{thm} \label{thm:icsine}
If $f$ is a function on a connected subset $E$ of $\Acal$, satisfying $f''(z)=-f(z)$ for all $z \in E$ and $f(0)=0$, $f'(0)=b \in\Acal$, then
$f(z)=b \star sin(z)$ for all $z \in E$.
\end{thm}

\noindent
{\bf Proof:} Let 
\begin{equation}
U(z)=f(z) \star sin(z)+f'(z) \star cos(z) \ \ \ \& \ \ \ V(z)=f(z) \star cos(z)-f'(z) \star sin(z).
\end{equation}
Apply the product rule to obtain:
\begin{align}
U'(z) &=f'(z) \star sin(z)+f(z) \star cos(z)+f''(z) \star cos(z)-f'(z) \star sin(z) \\ \notag
        &=f(z) \star cos(z)-f(z) \star cos(z)=0.
\end{align}
and
\begin{align} 
V'(z) &=f'(z) \star cos(z)-f(z) \star sin(z)-f''(z) \star sin(z)-f'(z) \star cos(z) \\ \notag
        &=-f(z) \star sin(z)+f(z) \star sin(z)=0.
\end{align}    
Thus, as $E$ is connected, $U$ and $V$ are constant.  Thus $U(z)=U(0)=b$ and $V(z)=V(0)=0$ for all z in E.  Hence,
\begin{equation}\label{eqn:forb}
b=f(z) \star sin(z)+f'(z) \star cos(z) 
\end{equation}
also
\begin{equation} \label{eqn:forcossin}
f(z) \star cos(z)=f'(z) \star sin(z).
\end{equation}
Combining Equations \ref{eqn:forcossin} and \ref{eqn:forb} with Theorem \ref{thm:pythagcosinesine} we derive:
\begin{align}
b \star sin(z) &=[f(z) \star sin(z)+f'(z) \star cos(z)] \star sin(z) \\ \notag
&=f(z)\star sin^2(z)+f'(z) \star sin(z) \star cos(z) \\ \notag
&=f(z)\star sin^2(z)+f(z)\star cos(z) \star cos(z) \\ \notag
&=f(z) \star [sin^2(z)+cos^2(z)] \\ \notag
&= f(z). \ \ \Box
\end{align}

\begin{thm} \label{thm:icsincos}
If $f$ is a function on a connected subset $E$ of $\Acal$, satisfying $f''(z)=-f(z)$ for all $z\in E$ and $f(0)=a$, $f'(0)=b\in\Acal$, then
$f(z)=a \star cos(z)+b \star sin(z)$ for all $z\in E$.
\end{thm}

\noindent
{\bf Proof:} Let $g(z)=f(z)-a\star cos(z)$. Then $g''(z)=-g(z)$ for all $z \in E$ and $g(0)=0$, $g'(0)=b$, so by Theorem \ref{thm:icsine}, we find $g(z)=b \star sin(z)$ thus $f(z)=a \star cos(z)+b\star sin(z). \, \square$ \\

\noindent
We now arrive at the angle addition formulas for trigonometric functions on $\Acal$:
\begin{thm}
For $z,w \in\Acal$,
$${\bf (i.)} \  \ sin(z+w)=sin(z)\star cos(w)+sin(w) \star cos(z)$$ 
$$ {\bf (ii.)} \ \ cos(z+w)=cos(z) \star cos(w)-sin(z)\star sin(w) $$
\end{thm}

\noindent
{\bf Proof:} Fix $w \in \Acal$ and let $f(z)=sin(z+w)$ for all $z \in\Acal$. Then $f''(z)=-f(z)$ for all $z \in\Acal$ and $f(0)=sin(w)$ and $f'(0)=cos(w)$, so by Theorem \ref{thm:icsincos}, we have
\begin{equation}
f(z)=sin(w) \star cos(z)+cos(w) \star sin(z) 
\end{equation}
for all $z,w \in\Acal$ which proves {\bf (i.)}. Continue to hold $w$ fixed and differentiate {\bf (i.)} with respect to $z$ to find: 
\begin{equation}
cos(z+w)=cos(z) \star cos(w)-sin(z) \star sin(w) 
\end{equation}
thus {\bf (ii.)} holds true. $\square$ \\

\noindent
These are helpful to uncover the component function content of sine or cosine over $\Acal$. 

\begin{exa}
Consider $\Hcal = \RN \oplus j \RN$ where $j^2=1$. For $x+jy \in \Hcal$ we calculate:
\begin{equation}
cos( x+jy) = cos(x) cos(jy) - sin(x) sin(jy).
\end{equation}
But, as $(jy)^{2n} = j^{2n}y^{2n} = y^2$ and $(jy)^{2n+1} = (j)^{2n+1} y^{2n+1} = jy^{2n+1}$ hence $\cos(jy) = \cos(y)$ 
and $ sin(jy) = j sin(y)$. Consequently,
\begin{equation}
cos( x+jy) = cos(x) cos(y) - j sin(x) sin(y).
\end{equation}
Differentiate with respect to $x$ holding $y$ fixed and find
\begin{equation}
sin( x+jy) = sin(x) cos(y) + j cos(x) sin(y).
\end{equation}
You might recognize the products of sine and cosine as well-known solutions to the unit-speed wave-equation $\phi_{xx}=\phi_{yy}$. This is no accident, the unit-speed wave equation is the generalized Laplace equation for $\Hcal$ and we know the component functions of an $\Hcal$-differentiable function solve the generalized Laplace equation of $\Hcal$. 
\end{exa}

\section{The N-Pythagorean theorem}
 \label{sec:nthagtheorem}
\noindent
In Section \ref{sec:transfunct} we studied functions whose properties were not tied to a particular choice of algebra. We saw how exponentials, cosine, sine, cosh and sinh all enjoy properties which hold in a multitude of algebras. The direction of the current section is quite the opposite. We now consider a method of obtaining new functions which are particular to our choice of algebra. We call such functions the {\it special functions} of $\Acal$.

\subsection{Special functions of an algebra}
The exponential function on an algebra is naturally defined by
\begin{equation}
e^z = 1+z+\frac{1}{2} z^2+ \cdots = \sum_{k=0}^{\infty} \frac{z^k}{k!}.
\end{equation}
We say the component functions of the exponential are the  {\em special functions} of the algebra. We explain how to calculate the special functions of a particular type of algebra in this section.
\begin{exa}
If $z \in \CN$ then $z = x+iy$ where $i^2=-1$. Notice the map $t \mapsto e^{it}$ is the composite of the complex exponential and the path $t \mapsto it$ in $\CN$. We calculate, 
\begin{align}
e^{it}=\sum_{n=0}^\infty \frac{(it)^n}{n!}
&=\sum_{n=0}^\infty i^{2n}\frac{t^{2n}}{(2n)!}+\sum_{n=0}^\infty i^{2n+1}\frac{t^{2n+1}}{(2n+1)!}\\ \notag
&=\sum_{n=0}^\infty (-1)^n\frac{t^{2n}}{(2n)!}+i\sum_{n=0}^\infty (-1)^n\frac{t^{2n+1}}{(2n+1)!}\\ \notag
&=cos(t)+isin(t).
\end{align}
The real and imaginary parts of this exponential are the real sine and cosine functions. Extending $t \in \RN$ to $z \in \CN$ provides the special functions $z \mapsto \cos(z)$ and $z \mapsto \sin(z)$ for $\CN$. These extended functions are the {\it special functions of $\CN$}
\end{exa}

\begin{exa}
Let $\mathcal{H}$ denote the hyperbolic numbers. If $z \in \mathcal{H}$ then $z = x+jy$ where $j^2=1$ and $x,y \in \RN$. Let $t \in \RN$ and calculate:
\begin{align}
e^{jt}=\sum_{n=0}^\infty \frac{(jt)^n}{n!}
&=\sum_{n=0}^\infty j^{2n}\frac{t^{2n}}{(2n)!}+\sum_{n=0}^\infty j^{2n+1}\frac{t^{2n+1}}{(2n+1)!}\\ \notag
&=\sum_{n=0}^\infty \frac{t^{2n}}{(2n)!}+j\sum_{n=0}^\infty \frac{t^{2n+1}}{(2n+1)!}\\ \notag
&=cosh(t)+jsinh(t).
\end{align}
Extending hyperbolic cosine and sine to $\Hcal$ we obtain the {\it special functions of $\Hcal$} are given by $z \mapsto \cosh(z)$ and $z \mapsto \sinh(z)$ for $z \in \Hcal$. 
\end{exa}

\noindent
Let us generalize the observations above for a unital algebra with generator $\varepsilon$. If each $z \in \Acal$ has the form $z = x_1+x_2\varepsilon+ \cdots + x_N \varepsilon^{N-1}$ for $x_1, x_2, \dots , x_N \in \RN$ then we say $\Acal$ is generated by $\varepsilon$.

\begin{thm} \label{thm:specialfunctions}
If $\Acal$ is an $N$-dimensional real algebra generated by $\varepsilon$ then there exist unique functions $f_1,\ f_2,\dots ,f_N: \RN \rightarrow \RN$ for which 
$$e^{\varepsilon t}=f_1(t)+\varepsilon f_2(t)+\dots +\varepsilon^{N-1}f_N(t) $$
for all $t \in \RN$. Moreover, for each $i=1,\dots , n$, there exist real constants $c_{ij}$ such that $f_i(t) = \sum_{j=0}^{\infty} c_{ij}t^j$ for all $t \in \RN$. That is, $f_1, \dots , f_n$ are entire on $\RN$. Furthermore, the functions $f_1, \dots , f_n$ uniquely extend to $\Acal$ via 
$f_i(z) = \sum_{j=0}^{\infty} c_{ij}z^j$ for each $z \in \Acal$.
\end{thm}

\noindent
{\bf Proof:} the map $z \mapsto \text{exp}(z)$ is an entire function on $\Acal$. It follows that $f(t) = exp(\varepsilon t) = \sum_{k=0}^{\infty} \frac{(\varepsilon t)^k}{k!}$ converges for each $t \in \RN$. We define components of $f$ by:
\begin{equation}
f(t) = exp( \varepsilon t) = f_1(t)+ \varepsilon f_2(t)+ \cdots + \varepsilon^{N-1} f_{N}(t) 
\end{equation}
for each $t \in \RN$. Similarly, we let $f_i^m: \RN \rightarrow \RN$ be the component functions of the map $t \mapsto \sum_{k=0}^{m} \frac{(\varepsilon t)^k}{k!}$. In particular,
\begin{equation}
\sum_{k=0}^m \frac{(\varepsilon t)^k}{k!} = f_1^m(t)+ \varepsilon f_2^m(t)+ \cdots + \varepsilon^{N-1} f_{N}^m(t). 
\end{equation}
Since $\sum_{k=0}^m \frac{(\varepsilon t)^k}{k!}  \rightarrow exp( \varepsilon t)$ as $m \rightarrow \infty$ we find $f_i^m(t) \rightarrow f_i(t)$ as $m \rightarrow \infty$ for each $t \in \RN$. It follows $f_1, f_2, \dots , f_N$ are entire on $\RN$; there exist real constants $c_{ij} \in \RN$ for which $f_i(t) = \sum_{j=0}^{\infty} c_{ij}t^j$ for all $t \in \RN$.  We conclude by applying Theorem \ref{thm:entire} which provides $f_i$ extends uniquely to an entire function on $\Acal$ for $i=1,2,\dots , N$. $\Box$

\begin{de}
Given $\Acal$ and $f_1, \dots , f_N: \Acal \rightarrow \Acal$ as discussed in Theorem \ref{thm:specialfunctions} we say $f_1, f_2, \dots , f_N$ are the {\em special functions} of $\Acal$.
\end{de}

\begin{exa}
If $\Hcal_3$ has basis $\{1,j,j^2\}$, where $j^3=1$.  We say $\Hcal_3$ are the $3$-hyperbolic numbers. Let $t \in \RN$, 
\begin{align}
e^{jt}&=\sum_{n=0}^\infty \frac{(jt)^n}{n!}\\ \notag
&=\sum_{n=0}^\infty j^{3n}\frac{t^{3n}}{(3n)!}+\sum_{n=0}^\infty j^{3n+1}\frac{t^{3n+1}}{(3n+1)!}+\sum_{n=0}^\infty j^{3n+2}\frac{t^{3n+2}}{(3n+2)!}\\ \notag
&=\sum_{n=0}^\infty \frac{t^{3n}}{(3n)!}+j\sum_{n=0}^\infty\frac{t^{3n+1}}{(3n+1)!}+j^2\sum_{n=0}^\infty\frac{t^{3n+2}}{(3n+2)!}.
\end{align}
Therefore, the special functions of the $3$-hyperbolic numbers are defined by:
\begin{equation} \label{eqn:larrymoeandcurly}
f_1(z)=\sum_{n=0}^\infty \frac{z^{3n}}{(3n)!}, \ \ 
f_2(z)=\sum_{n=0}^\infty\frac{z^{3n+1}}{(3n+1)!}, \ \ 
f_3(z)=\sum_{n=0}^\infty\frac{z^{3n+2}}{(3n+2)!}.
\end{equation}
\end{exa}

\begin{exa}
Suppose $\Gamma_3$ is generated by $\varepsilon$ with $\varepsilon^3=0$. The series for the exponential truncates nicely:
\begin{equation}
e^{\varepsilon x}=\sum_{n=0}^\infty \frac{(\varepsilon x)^n}{n!}=1+\varepsilon x+\frac{1}{2}\varepsilon^2 x^2.
\end{equation}
Thus we find special functions for the $3$-null numbers are
\begin{equation}
f_1(z)=1,\ \ \ \ f_2(z)=z, \ \ \ \ f_3(z) = \frac{1}{2}z^2 
\end{equation}
for all $z \in \RN \oplus \varepsilon \RN \oplus \varepsilon^2 \RN$.
\end{exa}

In similar fashion,the $\Gamma_N$ generated by $\varepsilon$ with $\varepsilon^N=0$ produces monomials $\frac{1}{k!}z^k$ for $k=0,1,\dots, N-1$. 

\section{N-trigonometric and N-hyperbolic functions}
Our goal in this section is to describe the special functions of $\CN_N$ and $\Hcal_N$. If $z \in \CN_N$ then we say $z$ is an $N$-complex number and
\begin{equation}
z = x_1+x_2j+\cdots + x_N j^{N-1}
\end{equation}
where $j^N=-1$ and $x_1, \dots , x_N \in \RN$. We also define $\Hcal_N$ to be the set of $N$-hyperbolic numbers which have the form
\begin{equation}
z = x_1+x_2j+\cdots + x_N j^{N-1}
\end{equation}
where $j^N=1$ and $x_1, \dots , x_N \in \RN$. 
Observe, 2-complex numbers are the ordinary complex numbers and 2-hyperbolic numbers form the hyperbolic numbers. 

\begin{de}
For a given $N\in\NN$, the $N$-trigonometric functions are defined on an algebra $\Acal$ as follows: $$  cos_N(z)=\sum_{k=0}^\infty (-1)^k\frac{z^{Nk}}{(Nk)!}  \ \ \& \ \
 sin_{N,p}(z)=\sum_{k=0}^\infty (-1)^k\frac{z^{Nk+p}}{(Nk+p)!} $$
for $p=1,2,\dots , N-1$. Likewise, the $N$-hyperbolic functions are defined by:
$$ cosh_N(z)=\sum_{k=0}^\infty\frac{z^{Nk}}{(Nk)!} \ \ \& \ \  sinh_{N,p}(z)=\sum_{k=0}^\infty\frac{z^{Nk+p}}{(Nk+p)!} $$
 for $p=1,2,\dots , N-1$.
\end{de}

\noindent
The series which define the $N$-trigonometric and $N$-hyperbolic functions converge for each $z \in \Acal$ for any algebra. For example, in complex analysis the hyperbolic functions are {\em entire} functions which have interesting applications. We pause to note the differential relations,
\begin{align} 
&\frac{d}{dz} \cos_N(z) = -\sin_{N,N-1}(z), 
&\frac{d}{dz} \cosh_N(z) = \sinh_{N,N-1}(z),\\ \notag
  &\frac{d}{dz} \sin_{N,1}(z) = \cos_{N}(z),
    &\frac{d}{dz} \sinh_{N,1}(z) = \cosh_{N}(z), \\ \notag
  &\frac{d}{dz} \sin_{N,p}(z) = \sin_{N,p-1}(z),
    &\frac{d}{dz} \sinh_{N,p}(z) = \sinh_{N,p-1}(z). 
\end{align}
for $p=2, \dots , N-1$ hold independent of our choice of $\Acal$. 

\begin{thm} We observe that:
\begin{enumerate}[{\bf (i.)}]
\item the N-trigonometric functions are the special functions of $\CN_N$, 
\item the N-hyperbolic functions are the special functions of $\Hcal_N$.
\end{enumerate}
\end{thm}

\noindent
{\bf Proof:} Begin with {\bf (i.)}. Suppose $j^N=-1$.  We have, for all $t \in \RN$:
\begin{align} \notag
e^{jt}
&=\sum_{k=0}^\infty \frac{(jt)^{Nk}}{(Nk)!}+\sum_{k=0}^\infty \frac{(jt)^{Nk+1}}{(Nk+1)!}+\dots +\sum_{k=0}^\infty \frac{(jt)^{Nk+N-1}}{(Nk+N-1)!}\\ \notag
&=\sum_{k=0}^\infty \frac{(j^N)^kt^{Nk}}{(Nk)!}+\sum_{k=0}^\infty \frac{j(j^N)^kt^{Nk+1}}{(Nk+1)!}+\dots +\sum_{k=0}^\infty \frac{j^{N-1}(j^N)^kt^{Nk+N-1}}{(Nk+N-1)!}\\ \notag
&=\sum_{k=0}^\infty \frac{(-1)^kt^{Nk}}{(Nk)!}+j\sum_{k=0}^\infty \frac{(-1)^kt^{Nk+1}}{(Nk+1)!}+\dots +j^{N-1}\sum_{k=0}^\infty \frac{(-1)^kt^{Nk+N-1}}{(Nk+N-1)!}\\ \notag
&=\cos_N(t)+j\sin_{N,1}(t)+\dots+j^{N-1}\sin_{N,N-1}(t).
\end{align}
Thus $\cos_N,\ \sin_{N,1},\dots, \sin_{N,N-1}$ form the special functions of $\CN_N$. In contrast, for $\Hcal_N$ we have $j^N=1$ thus $j^{Nk+p} = (j^N)^kj^p = j^p$ and it follows we find the special functions of $\Hcal_N$ are $\cosh_N,\ \sinh_{N,1},\dots, \sinh_{N,N-1}$. $\Box$

\section{Pythagorean functions}
We may study a finite dimensional associative algebras over $\RN$ by instead studying the associated {\em regular representation} of $\Acal$. When $\Acal = \RN^N$ the matrix corresponding to $z \in \Acal$ is the standard matrix of the left multilication by $z$ map; $L_z: \Acal \rightarrow \Acal$ is defined by $L_z(x) = zx$ for all $x \in \Acal$ and 
\begin{equation}
 \mathbf{M}(z) = [L_z] = \left[ ze_1 | ze_2| \cdots | ze_N \right]. 
\end{equation}
When we use $e_1=1 \in \Acal$ the matrix representation of $z$ is further simplified.

\begin{exa}
If $z = x+iy \in \CN$ then note $e_1=1$ and $e_2=i$ hence
\begin{equation}  
\mathbf{M}(z) = [z|zi] = [x+iy|xi-y] = \left[ \begin{array}{cc} x & -y \\ y & x \end{array}\right] 
\end{equation}
\end{exa}

\begin{exa}
If $z \in \Hcal_N$ where $e_1=1, e_2=j, \dots , e_N = j^{N-1}$ and $j^N=1$ then
\begin{equation}  
\mathbf{M}(z) = [z|zj| \cdots | zj^{N-1}] 
\end{equation}
To be explicit, in the $N=3$ case we have:
\begin{equation}  
\mathbf{M}(x+jy+zj^2) = \left[ \begin{array}{ccc} x & z & y \\
y & x & z \\
z & y & x\end{array}\right]
\end{equation}
\end{exa}



Combining the exponential on $\Acal$ generated by $j$ with the determinant function allows us to create a new function on $\Acal$ which we define below:

\begin{de}
For $\Acal = \RN^N$ with basis $1, j , j^2, \dots , j^{N-1}$ we define the {\em real Pythagorean Function} by
$$ \mathcal{P}_{\Acal}(t) = \text{det}\left[ e^{j t}| j e^{j t} | \cdots | j^{N-1}  e^{j t} \right]$$
for each $t \in \RN$.  The {\em Pythagorean function of $\Acal$} is the unique extension of the real Pythagorean function to $\Acal$. \end{de}

\noindent
Notice Pythagorean function is manifestly a formula which involves the special functions of $\Acal$. If one is willing to study the matrices with components taken from an algebra then we may express
$$ \mathcal{P}_{\Acal}(z) = \text{det}\left[ e^{j z}|j  e^{j z} | \cdots | j^{N-1}  e^{j z} \right]$$
for each $z \in \Acal$. Once more it is instructive to examine how this construction unfolds for complex and hyperbolic numbers. 

\begin{exa}
In $\CN$, as $e^{it} = \cos( t) + i \sin( t)$ and $ie^{it} = i \cos( t)-\sin( t)$ we find the Pythagorean function:
\begin{equation}
\mathcal{P}_{\CN}(z) =det\left[ \begin{array}{rr} \cos (z)& -\sin (z) \\ \sin( z) & \cos (z)\end{array}\right]
= \cos^2 (z) +\sin^2 (z).
\end{equation}
Observe $\mathcal{P}_{\CN}(z)  = 1$ for all $z \in \CN$. 
\end{exa}

\begin{exa}
In $\Hcal$ which is generated by $j$ with $j^2=1$ we found $e^{jt} = \cosh(t) + j \sinh(t)$ and $je^{jt} = j\cosh(t)+\sinh(t)$. Thus, the Pythagorean function is:
\begin{equation}
\mathcal{P}_{\Hcal}(z) =det\left[ \begin{array}{rr} cosh(z)& sinh(z)\\ sinh(z) & cosh(z)\end{array}\right]
=cosh^2(z)-sinh^2(z). 
\end{equation}
It can be shown $\mathcal{P}_{\Hcal}(z)  = 1$ for all $z \in \Hcal$. 
\end{exa}

\noindent
These examples are part of a larger pattern: 

\begin{thm}
\textbf{(N-Pythagorean Theorem)} If $\Acal=\RN^N$ is generated by $j$ where $j^N = c \in \RN$ then $\mathcal{P}_{\Acal}(z) = 1$ for all $z \in \Acal$.
\end{thm}

\noindent
We prove this result in Section \ref{sec:proofofNpythag}. Our examples are based on the choices $c = \pm 1$ or $c=0$. However, the choice $c=0$ is not especially interesting:

\begin{exa}
Consider the dual numbers $\Acal = \RN \oplus \eta \RN$ with $\eta^2=0$. Observe, $e^{\eta t} = 1+\eta t$ hence $\eta e^{\eta t} = \eta$ and so 
\begin{equation}
\mathcal{P}_{\Acal}(z) = \text{det} \left[ \begin{array}{cc} 1 & 0 \\ z & 1 \end{array}\right] = 1.
\end{equation}
\end{exa}

In contrast, the identities found below are surprising.

\begin{exa}
In $\Hcal_3$ we have basis $1,j,j^2$ with $j^3=1$ and we calculate:
\begin{align}
\mathcal{P}_{\Hcal_3}(z) &=\text{det} \left[ \begin{array}{ccc} cosh_3(z) & sinh_{3,2}(z)& sinh_{3,1}(z)\\ sinh_{3,1}(z)& cosh_3(z)& sinh_{3,2}(z)\\ sinh_{3,2}(z)& sinh_{3,1}(z)& cosh_3(z)\end{array} \right]\\ \notag
&=cosh_3^3(z)+sinh_{3,1}^3(z)+sinh_{3,2}^3(z)-3cosh_3(z)sinh_{3,1}(z)sinh_{3,2}(z).
\end{align}
The $N$-Pythagorean Theorem indicates $\mathcal{P}_{\Hcal_3}(z) = 1$ for each $z \in \Hcal_3$.
\end{exa}

\begin{exa}
Consider $\CN_3$ generated by $j$ with $j^3=-1$. Then $e^{jz} = \cos_3(z)+j\sin_{3,1}(z)+ j^2 \sin_{3,2}(z)$ and we calculate:
\begin{align}
\mathcal{P}_{\CN_3}(z) &=\text{det} \left[ \begin{array}{ccc} 
\cos_3(z)		& -\sin_{3,2}(z)	& -\sin_{3,1}(z)	\\ 
\sin_{3,1}(z)	& \cos_3(z)		& -\sin_{3,2}(z)	\\
\sin_{3,2}(z)	& \sin_{3,1}(z)	& \cos_3(z)		
\end{array} \right]\\ \notag
&=cos_3^3(z)-sin_{3,1}^3(z)+sin_{3,2}^3(z)+3cos_3(z)sin_{3,1}(z)sin_{3,2}(z).
\end{align}
The $N$-Pythagorean Theorem indicates $\mathcal{P}_{\CN_3}(z) = 1$ for each $z \in \CN_3$.
\end{exa}

The relations implied by $\mathcal{P}_{\Acal}(z)=1$ for $N \geq 4$ are rather lengthy. 

\section{Proof of N-Pythagorean theorem}
\label{sec:proofofNpythag}

\noindent
{\bf Proof:} Assume $\Acal = \RN^N$ is generated by $j$ with $j^N=c$ for some given $c \in \RN$. Begin by writing an explicit formula for the real Pythagorean function:
\begin{equation}
 \mathcal{P}_{\Acal}(t) = \text{det} ( \mathbf{M}(e^{jt}) )= \text{det} \left[ e^{jt}|je^{jt}| \cdots | j^{N-1}e^{jt} \right]. 
\end{equation}
If $f: \text{dom}(f) \subseteq \Acal \rightarrow \Acal$ is $\Acal$-differentiable and $g: \RN \rightarrow \Acal$ is real differentiable then there is a chain rule for the composite $f \comp g: \RN \rightarrow \Acal$. In particular, see Theorem 6.7 of~\cite{cookAcalculusI},
\begin{equation}
 \frac{d}{dt} (f (g(t)) = \frac{df}{dz}(g(t)) \frac{dg}{dt}.
\end{equation}
Here $g = g_1+ jg_2+ \cdots + j^{N-1}g_N$ has $\frac{dg}{dt} = \frac{dg_1}{dt}+ j\frac{dg_2}{dt}+ \cdots + j^{N-1}\frac{dg_{N-1}}{dt}$. Consider $f(z) = e^z$ and $g(t) = jt$ then we know $\frac{df}{dz} = e^z$ and it is simple to calculate $\frac{dg}{dt} = j$. Hence,
\begin{equation}
 \frac{d}{dt} (e^{jt}) = je^{jt}. 
\end{equation}
Multiply by $j^{p-1}$ and note:
\begin{equation}
 j^{p-1}\frac{d}{dt} (e^{jt}) = j^{p-1}je^{jt} \ \ \Rightarrow \ \ \frac{d}{dt} \left(j^{p-1}e^{jt} \right) = j^pe^{jt}. 
\end{equation}
In particular, $\frac{d}{dt} \text{Col}_{p}(\mathbf{M}(e^{jt}) ) = \text{Col}_{p+1}(\mathbf{M}(e^{jt}) )$ for $p=1, \dots , N-1$. However,
\begin{equation}
 \frac{d}{dt} (j^{N-1}e^{jt}) = j^Ne^{jt}  
\end{equation}
provides $\frac{d}{dt} \text{Col}_{N}(\mathbf{M}(e^{jt}) ) = c\, \text{Col}_{1}(\mathbf{M}(e^{jt}) )$ where $j^N = c$. Hence observe every column in the regular representation of $e^{jt}$ has a derivative which is proportional to another column in the representation. \\

\noindent
In order to compress the notation a bit let us set
$\mathbf{M}(e^{jt}) = B = [B_1|\cdots |B_N]$ which gives $\mathcal{P}_{\Acal}(t) = \text{det}(B)$. The formula for the determinant below makes manifest the fact the determinant is a multilinear function of its columns:
\begin{equation}
 \text{det}(B) = \sum_{i_1\dots i_N=1}^{N} \epsilon_{i_1i_2\dots i_N}B_{i_11}B_{i_22} \cdots B_{i_NN} 
\end{equation}
here $\epsilon_{i_1\dots i_N}$ is the completely antisymmetric symbol where $\epsilon_{12\dots N}=1$. By the $N$-fold product rule we find: $\frac{d}{dt}( \text{det}(B)) =$
\begin{align} 
   &= \sum_{i_1\dots i_N} \epsilon_{i_1\dots i_N}\frac{dB_{i_11}}{dt}B_{i_22} \cdots B_{i_NN}+ \sum_{i_1\dots i_N} \epsilon_{i_1\dots i_N} B_{i_11}\frac{dB_{i_22}}{dt} \cdots B_{i_NN} \\ \notag
& \qquad + \cdots +  \sum_{i_1\dots i_N} \epsilon_{i_1\dots i_N} B_{i_11}B_{i_22} \cdots \frac{dB_{i_NN}}{dt} \\ \notag
&= \text{det}[B_2|B_2|\cdots |B_N]+
\text{det}[B_1|B_3|B_3|\cdots |B_N]+ \cdots + \text{det}[B_1|B_2|\cdots |c B_1] \\ \notag
&=0.
\end{align}
Thus $t \mapsto \text{det} ( \mathbf{M}(e^{jt}) )$ is a constant function on $\RN$. Notice $t=0$ maps to $\text{det}(I)=1$ hence $\mathcal{P}_{\Acal}(t) = 1$. for each $t \in \RN$.  Finally, we extend to all of $\Acal$ using Theorem \ref{thm:entire}. $\Box$

\section{Acknowledgements}
The authors are thankful to N. BeDell for helpful comments on a rough draft of this article. We should mention that W.S. Leslie provided an alternate, purely algebraic, proof of the $N$-Pythagorean Theorem in private communication. Finally, we thank Khang Nguyen for finding the improved Theorem 5.4 as well as reformulations of Corollary 5.20 and 5.21.

\end{document}